\documentclass[11pt]{article}

\usepackage[german,english]{babel}
\usepackage[cp850]{inputenc}
\usepackage{latexsym,graphicx}
\usepackage{amsmath}
\usepackage{color}

\font\tenmath=msbm10 scaled 1200
\font\sevenmath=msbm7 scaled 1200
\font\fivemath=msbm5 scaled 1200

\newtheorem{Lemma}{Lemma}[section]

\newfam\mathfam \textfont\mathfam=\tenmath
\scriptfont\mathfam=\sevenmath \scriptscriptfont\mathfam=\fivemath
\def\math{\fam\mathfam}
\def\R{{\math R}}
\def\N{{\math N}}
\def\E{{\math E}}

\def\P{{\math P}}

\def\Q{{\math Q}}
\def\B{{\math B}}

\def \F{{\cal F}}

\def \^#1{\if#1i{\accent"5E\i}\else{\accent"5E#1}\fi}

\def \cqfd{\quad_\diamondsuit}
\def \ms{\medskip}
\def \ss{\smallskip}
\def \bs{\bigskip}
\def \ni{\noindent}

\def \supess{{\rm esssup}}

\newtheorem{Thm}{Theorem}

\newtheorem{Pro}{Proposition}
\newtheorem{Cor}{Corollary}

\author{
{\sc Olivier Bardou}\thanks{Gaz de France, Research and Development
Division, 361 Avenue du Pr\'esident Wilson - B.P. 33, 93211
Saint-Denis La Plaine cedex. E-mail: {\tt
olivier-aj.bardou@gazdefrance.com}} \quad{\sc Sandrine Bouthemy}\thanks{Gaz de France, Research and Development Division,
361 Avenue du Pr\'esident Wilson - B.P. 33, 93211 Saint-Denis La
Plaine cedex. E-mail: {\tt sandrine.bouthemy@gazdefrance.com}} \quad {\sc  and}
\quad {\sc Gilles Pag\`es} \thanks{Laboratoire de Probabilit\'es et Mod\`eles al\'eatoires, UMR~7599, Universit\'e Paris 6, case 188, 4,
pl. Jussieu, F-75252 Paris Cedex 5, France. E-mail:{\tt  gpa@ccr.jussieu.fr}}
}
\date{}
\title{\bf When are Swing options bang-bang and how to use it?}
\date{6th April 2007}
\begin{document}

\maketitle
\begin{abstract}
In this paper we investigate a class of swing options with firm
constraints in view of the modeling of supply agreements.  We show,
for a fully general payoff process, that the premium, solution to a
stochastic control problem, is concave and piecewise affine as a
function of the global constraints of the contract. The existence of
bang-bang optimal controls is established for a set of constraints
which generates by affinity the whole premium function.

When the payoff process is driven by an underlying Markov process,
we propose a quantization based recursive backward procedure to
price these contracts. A priori error bounds are established,
uniformly with respect to the global constraints.
\end{abstract}


\bigskip
%
\noindent
\section{Introduction}

\setcounter{equation}{0}
\setcounter{Assumption}{0}
\setcounter{Theorem}{0}
\setcounter{Proposition}{0}
\setcounter{Corollary}{0}
\setcounter{Lemma}{0}
\setcounter{Definition}{0}
\setcounter{Remark}{0}
The deregulation of energy markets has given rise to various
families of contracts. Many of them appear as some derivative
products whose underlying is some tradable futures (day-ahead, etc)
on gas or electricity (see~\cite{GEM} for an introduction). The
class of swing options has been paid a special attention in the
literature, because it includes many of these derivative products. A
common feature to all these options is that they introduce some risk
sharing between a producer and a trader, of gas or electricity for
example. From a probabilistic viewpoint, they appear as some
stochastic control problems modeling multiple optimal stopping
problems (the control variable is the purchased quantity of energy);
see $e.g.$~\cite{CATO, CADA} in a continuous time setting. Gas
storage contracts (see~\cite{GobetDF},~\cite{CALU1}) or electricity
supply agreements (see~\cite{KEP}, \cite{CALU2}) are examples  of
such swing options. Indeed, energy supply contracts are one simple
and important example of such swing options that will be deeply
investigated in this paper (see below, see also~\cite{GEM} for an
introduction). It is worth mentioning that this kind of contracts
are slightly different from multiple exercises American options as
considered in~\cite{CATO} for example. In our setting the volumetric
constraints play a key role and thus, the flexibility is not
restricted to time decisions, but also has to take into account
volumes management.

Designing efficient numerical procedures for the pricing of swing
option contracts remains a very challenging question as can be
expected from a possibly multi-dimensional stochastic control
problem subject to various constraints (due to the physical
properties of the assets like in storage contracts). Most recent
approaches developed in mathematical finance, especially for the
pricing of American options, have been adapted and transposed to the
swing framework: tree (or ``forest") algorithms in the pioneering
work~\cite{JAROTO}, Least squares regression MC methods
(see~\cite{GobetDF}), PDE's numerical methods (finite elements,
see~\cite{WIWI}).

The aim of this paper is to deeply investigate an old question,
namely to elucidate the structure of the optimal control in supply
contracts (with firm constraints) and how it impacts the numerical
methods of pricing. We will provide  in  a quite general (and
abstract) setting some ``natural" (and simple) conditions involving
the local and global purchased volume constraints  to ensure the
existence of {\em bang-bang optimal strategy} (such controls usually
do not exist). It is possible to design  {\em a priori} the contract
so that their parameters  satisfy these conditions. To our knowledge
very few theoretical results have been established so far on this
problem (see however~\cite{GobetDF} in a Markovian framework for
contracts with penalized constraints and~\cite{ROZH}, also in a
Markovian framework).

This first result of the paper not only enlightens the understanding
of the management of a swing contract: it also has some deep
repercussions on the numerical methods to price it. As a matter of
fact, taking advantage of the existence of a bang-bang optimal
strategy, we propose and analyze in details (when the underlying
asset has a Markovian dynamics) a quantized Dynamic Programming
procedure to price any swing options  whose volume constraints
satisfy the ``bang-bang" assumption. Furthermore some {\em a priori}
error bounds are established. This procedure turns out to be
dramatically efficient, as emphasized in the companion
paper~\cite{BABOPA1} where the method is extensively tested with
assets having multi-factor Gaussian underlying dynamics  and
compared to the least squares regression method.

\paragraph{The abstract swing contract with firm constraints}
The holder of a supply contract has the right to purchase
periodically (daily, monthly, etc) an amount of energy at a given
unitary price. This amount of energy is subject to some lower and
upper ``local" constraints. The total amount of energy purchased at
the end of the contract is also subject  to a ``global" constraint.
Given  dynamics on the energy price process, the problem is to
evaluate the price of such a contract, at time $t=0$ when it is
emitted and during its whole life up to its maturity.

To be precise, the owner of the contract is allowed to purchase at times $t_i$, $i=0,\ldots,n-1$ a quantity $q_{i}$ of energy at a
unitary {\em strike price} $K_i:=K(t_i)$. At every date $t_i$, the purchased quantity $q_{i}$ is  subject to the firm {\em ``local"
constraint},
\[
q_{\min}\le q_{i}\le q_{\max},\qquad i=0,\ldots,n-1,
\]
whereas the global purchased quantity $\bar q_{n}:= \sum_{i=0}^{n-1}q_{i}$ is subject to the (firm) global constraint
\[
\bar q_{n}\!\in [Q_{\min},Q_{\max}]\qquad (0<Q_{\min} \le Q_{\max}<+\infty).
\]

The strike price process $(K_i)_{0\le i\le n-1}$  can be either
deterministic (even constant) or stochastic, $e.g.$ indexed on
past values of other commodities (oil, etc). Usually, on energy
markets the price is known through future contracts
$(F_{s,t})_{0\le s\le t}$ where $F_{s,t}$ denotes the price at
time $s$ of the forward contract delivered at maturity $t$. The
available data at time $0$ are $(F_{0,t})_{0\le t\le T}$ (in real
markets this is of course not a continuum).

The underlying asset price process, temporarily denoted $(S_{t_i})_{0\le i\le n-1}$,  is often the so-called ``day-ahead" contract
$F_{t,t+1}$ which is a tradable instrument or the spot price $F_{t,t}$ which is not. All the decisions about the contract need to be
adapted to the filtration of $(S_{t_i})$ $i.e.$ ${\cal F}_i:=\sigma(S_{t_j}, j=0,\ldots,i)$, $i=0,\ldots,n-1$ (with ${\cal F}_0=\{\emptyset,\Omega\}$). This means that the price of such a contract is given at any time $t_k$, by
\begin{eqnarray*}
P^n_k(Q^k_{\min},Q^k_{\max})&:=& \supess\left\{\E\left(\sum_{j=k}^{n-1} \!q_{j}e^{-r(t_j-t_k)}(S_{t_j}-K_j)\,|\,{\cal
F}_k\right),\right.\\ &&\left. q_j:(\Omega,{\cal F}_j)\to [q_{\min}, q_{\max}],j=k,\ldots,n-1, \sum_{j=k}^{n-1}q_j\!\in [Q^k_{\min},
Q^k_{\max}]\right\}
\end{eqnarray*}
where $Q^k_{\min}= Q_{\min}-\bar q_k$, $Q^k_{\max}= Q_{\max}-\bar q_k$ denote the {\em residual global constraints} and $r$ denotes the
(deterministic) interest rate. This pricing problem clearly appears as a stochastic control problem.

In the pioneering work by~\cite{JAROTO}, this type of contract was computed by using some forests of (multinomial) trees.
 A natural variant, at least for numerical purpose, is to consider a penalized version of this stochastic control. Thus, in~\cite{GobetDF},
a penalization $Q_{n}(V_n, \bar q_{n})$ with $Q_{n} (x,q)= -((x-Q_{\max})^++(x-Q_{\min})^+)/\varepsilon$ is added ($Q_n$ is negative outside
$[Q_{\min},Q_{\max}]$ and zero inside).

As concerns more sophisticated contracts (like storages),  the
holder of
 the contract receive a quantity $\Psi(t_i, S_{t_i}, q_{i})$ when deciding $q_i$. When dealing with gas this is due to the storing
constraints since injecting or withdrawing gas from its storing units induce fixed costs (and physical constraints (pressure, etc)).

As concerns the underlying asset dynamics, it is commonly shared in finance to assume that the traded asset has a Markovian dynamics (or
is a component of a Markov process like with stochastic volatility models).  The dynamics of physical assets for many reasons (some of
them simply coming from history) are often modeled using some more deeply non-Markovian models like long memory processes, etc.

All these specific features of energy derivatives suggest to tackle
the above pricing problem in a rather general framework, trying to
avoid as long as possible to call upon Markov properties. This is
what we do in the first part of the paper where the general setting
of a   swing option defined by an abstract sequence of ${\cal
F}_{k}$-adapted payoffs is deeply investigated as a function of its
global constraints $(Q_{\min},Q_{\max})$ (when the local constraints
are normalized $i.e.$ $q_i$ is $[0,1]$-valued for every $i\!\in
\{0,\ldots,n-1\}$). We show that this premium is a concave,
piecewise affine, function  of the global constraints, affine on
triangles of the $(m,M)+\{(u,v),\; 0\le u\le v\le 1\}$, $m,\, M\!\in
\N^2$, $m\le M\le n$ and $(m,M)+\{(u,v),\; 0\le v\le u\le 1\}$,
$m,\, M\!\in \N^2$, $m\le M-1\le n-1$. We also show that for
integral valued global constraints, the optimal controls  are always
bang-bang $i.e.$ the {\em a priori} $[0,1]$-valued optimal purchased
quantities $q^*_i$ are in fact always equal to $0$ or $1$. Such a
result can be extended in some way to any couple of global
constraints when all the payoffs are nonnegative.

Then, when there is an underlying Markov ``structure process",   we propose an optimal quantization
based on numerical approach to price  efficiently swing options.  This Markov ``structure process" can
be the underlying traded asset
itself or a  higher dimensional  hidden  Markov process: such a framework comes out in case of multi-factor processes having some
long-memory properties.

Optimal Quantization was first introduced as a numerical method to solve nonlinear problem arising in Mathematical Finance in a series
of papers~\cite{BAPA0, BAPA1, BAPA2, BAPAPR} devoted to the pricing and hedging of American style multi-asset options.  It has also been applied to stochastic control problem, namely portfolio optimization in~\cite{PAPHPR}. The
purely numerical aspects induced by optimal quantization, with a special emphasis on the Gaussian distribution,  have been investigated in~\cite{PAPR}. See~\cite{PAPHPR2} for a survey on numerical application of optimal
quantization to Finance. For other applications (to automatic classification, clustering, etc), see also~\cite{GEGR}. In this paper, we propose
a quantized backward dynamic programming to  approximate the premium of a swing contract. We analyze the rate of convergence of this
algorithm and provide some {\em a priori} error bounds in terms of quantization errors.

We illustrate the method by computing the whole graph of the premium viewed as a function of the global constraints, combining the affine
property of the premium and the quantized algorithm in ``toy model":   the future prices of  gas  are modeled by a two factor Gaussian
model.   An extensive study of the pricing method by optimal
quantization is carried out from both a theoretical and numerical point of view in~\cite{BABOPA1}.

\bs
The paper is organized as follows. In the section below we detail the decomposition of swing options  into a swap contract and a normalized swing option. In
Section~\ref{abstractswing}, we precisely describe our abstract setting for normalized swing options with firm constraints and the variable of
interest (global constraints, local constraints, etc).  In Section~\ref{ProgDynAbst}, we establish  the dynamic  programming formula satisfied
in full generality by the premium as a function of the global constraints  (this unifies the similar results obtained in Markov settings,
see~\cite{JAROTO},~\cite{GobetDF}, etc) and  we show this is a  concave function with respect to the global constraints. Then, in
Section~\ref{piecaffine}, we prove   in our abstract framework  that the premium function is   piecewise affine   and that
the optimal purchased quantities satisfy a ``$0$-$1$" or bang-bang principle (Theorem~\ref{main}). A special
attention is paid to the
$2$-period model which provides an intuitive interpretation of the results. In Section~\ref{SwingQuant}, after some short background on
quantization and its optimization, we propose a quantization based backward dynamic programming formula as a numerical method to solve the
swing pricing problem. Then we provide some error bounds for the procedure depending on the quantization error induced by the quantization of
the Markov structure process.

\bs \ni{\sc Notations.} $\bullet$ The Lipschitz coefficient of a
function $f:\R^d\to \R$ is defined by  $\displaystyle [f]_{\rm
Lip}:=\sup_{x\neq y}\frac{|f(x)-f(y)|}{|x-y|}\le +\infty$.  The
coefficient $[f]_{\rm Lip}$ is finite if and only if $f$ Lipschitz
continuous.

$\bullet$ The canonical Euclidean norm on $\R^d$ will be denoted $|\,.\,|$.

\section{Canonical decomposition, normalized swing option}
As a first step we need to normalize this contract to reduce some useless technical aspects. In practice this normalization, in fact
decomposition, corresponds to the splitting of the contract into a swap and a normalized swing. The decomposition   can be derived from
the fact that an ${\cal A}$-measurable random variable $q$ is $[q_{\min},q_{\max}]$-valued  if and only if there exists a $[0,1]$-valued
${\cal A}$-measurable random variable $q'$ such that $ q= q_{\min}+ (q_{\max}-q_{\min})q'$.
Then, for every $k\!\in\{0,\ldots,n-1\}$,
\begin{eqnarray*}
P^n_k(Q^k_{\min},Q^k_{\max})&=&q_{\min}\underbrace{\sum_{j=k}^{n-1} \!e^{-r(t_j-t_k)}\E( S_{t_j}-K_j\,|\,{\cal
F}_k)}_{\mbox{\em
swap contract}} +(q_{\max}-q_{\min})\!\underbrace{P^{[0,1],n}_k\left(\widetilde Q^k_{\min},\widetilde Q^k_{\max}\right)}_{\mbox{\em
normalized contract}}
\end{eqnarray*}
where
$$
\widetilde Q^k_{\min}=\frac{Q^k_{\min}-(n-k)q_{\min}}{q_{\max}-q_{\min}} \qquad\mbox{ and }\qquad \widetilde
Q^k_{\max}=\frac{Q^k_{\max}-(n-k)q_{\min}}{q_{\max}-q_{\min}}
$$
and $P^{[0,1],n}_k\left(\widetilde Q^k_{\min},\widetilde Q^k_{\max}\right)$ is a {\em normalized
swing contract} in which the local constraints are $[0,1]$-valued .

\section{An abstract model for swing options with firm constraints}\label{abstractswing}

\setcounter{equation}{0}
\setcounter{Assumption}{0}
\setcounter{Theorem}{0}
\setcounter{Proposition}{0}
\setcounter{Corollary}{0}
\setcounter{Lemma}{0}
\setcounter{Definition}{0}
\setcounter{Remark}{0}
%

One considers a sequence $(V_k)_{0\le k\le n-1}$ of integrable random variables defined on a
probability space   $(\Omega,{\cal A},\P)$. Let ${\cal F}_k^V:=\sigma(V_0,V_1,\ldots,V_k),\; k=0,\ldots,n-1$ denote its natural
filtration. For convenience we introduce a more general discrete time  filtration
${\cal F}:=({\cal F}_k)_{0\le k\le n-1}$ to which $V=(V_k)_{0\le k\le n-1}$ is adapted $i.e.$ satisfying ${\cal F}_k^V\subset {\cal
F}_k$,
$k=0,\ldots,n-1$ $i.e.$ such that the sequence $(V_k)_{0\le k\le n-1}$ is  ${\cal F}$-adapted.

We aim to solve the following abstract stochastic  control problem (with {\em maturity} $n$)
\begin{equation}\label{swing}
({\cal S})^{\cal F}_{n}\equiv
\supess\left\{\E\left(\sum_{k=0}^{n-1} \!q_{k}V_k\,|\, {\cal
F}_0\right),\; q_k:(\Omega,{\cal F}_k)\to [0,1],0\le k\le
\ldots,n-1, \sum_{k=0}^{n-1}q_k\!\in [Q_{\min}, Q_{\max}]\right\}
\end{equation}

\noindent where $Q_{\min}$ and $Q_{\max}$ are two non-negative ${\cal F}_0$-measurable random variables
satisfying
\begin{equation}\label{globconstr}
0\le Q_{\min}\le Q_{\max}\le n .
\end{equation}

(The inequality $Q_{\max}\le n$ induces no loss of generality: one
can always replace $Q_{\max}$ by $Q_{\max}\wedge n$ in~(\ref{swing})
since in any case  $q_0+\dots+q_n\le n$).

Note that no
assumption is made on the dynamics of the state process
$V=(V_k)_{0\le k\le n-1}$.

We need to introduce the following notations and terminology:

-- Throughout the paper, $\supess$ will always be taken with respect to the probability $\P$  so  $\P$ will be dropped from now on.

\ss
 -- A couple $Q:=(Q_{\min},Q_{\max})$  of non-negative ${\cal F}_k$-measurable random variables satisfying
$0\le Q_{\min}\le Q_{\max}\le n\!-\!k$ is called   a  couple of {\em global constraints} at time $k$.

\ss
-- An $ {\cal F}$-adapted sequence     $q=(q_k)_{0\le k\le n-1}$ of $[0,1]$-valued r.v. is called a   {\em locally
admissible control}. For any locally admissible control, one defines the cumulative purchase process by
\[
\bar q_0:=0,\quad \bar q_k:=q_0+\cdots+q_{k-1},\; k=1,\ldots, n.
\]

If $\bar q_{n}\!\in [Q_{\min}, Q_{\max}]$, $q$ is called an {\em $({\cal F},Q)$-admissible control}.

\ss
-- For every $k\!\in \{0,\ldots,n\}$ and every couple of global constraints $Q:=(Q_{\min},Q_{\max})$ at time $k$,  set
\begin{equation}\label{(2.3)}
P^n_k(Q,(V,\F)):= \supess\left\{\E\left(\sum_{\ell=k}^{n-1}
q_{\ell}V_\ell\,|\, {\cal F}_k\right),\; q\; ({\cal
F},Q)\mbox{-admissible control}\right\}
\end{equation}
so that  $P^n_0(Q,(V,\F))$ is the value function of the stochastic
control problem~(\ref{swing}) when the global constraints
$Q_{\min}$ and $Q_{\max}$ (at time $0$)
satisfy~(\ref{globconstr}). Note that the standard convention
$\supess(\o)=0$ yields $P^n_{n}\equiv 0$. To alleviate notations, the payoff process $V$
the filtration ${\cal F}$ will be often dropped in
$P^n_k(Q,(V,\F))$.

\bs
To be precise we will answer the   following questions:

\ms \noindent $\bullet$  Existence of an optimal control
$q^*=(q^*_k)_{0\le k\le n-1}$.

\ms \noindent $\bullet$ Regularity of the value function $Q\mapsto
P^n_k(Q,(V,\F))$.

\ms \noindent $\bullet$ Existence of a {\em bang-bang} optimal
control $(q^*_k)_{0\le k\le n-1}$ for certain values of the global
constraints $Q$ (namely when $Q$ has integral components)~?

By bang-bang we mean that $\P(d\omega)$-$a.s.$

\ss
-- all the {\em local} constraints on the $q_k(\omega)$ are saturated $i.e.$  for every
$k\!\in\{0,\ldots,n-1\}$,
$q_k(\omega)\!\in\{0,1\wedge Q_{\max}\}$

or

\ss
-- there exists at most one instant $k_0(\omega)$  such that $q_{k_0(\omega)}(\omega)\!\in (0,1\wedge Q_{\max})$ and one  global
constraint is saturated.

\ss Note that if $Q\!\in \N^2$, then a bang-bang $Q$-admissible
control necessary satisfies $\P(d\omega)$-$a.s.$ $q_k(\omega)\!\in
\{0,1\}$, $k=0,\ldots, n-1$. The existence of bang-bang optimal
controls combined with the piecewise affinity of $P^n_k(Q)$ will
be  the key  in the design of a numerical.

\ms \noindent  $\bullet$ When there is  an underlying structure
Markov process ($V_k=v_k(Y_k))$, we will show that the optimal
control turns out to be a function of $Y_k$ at every time $k$ as
well.

%

\section{Abstract dynamical programming principle and first properties}\label{ProgDynAbst}
\setcounter{equation}{0}
\setcounter{Assumption}{0}
\setcounter{Theorem}{0}
\setcounter{Proposition}{0}
\setcounter{Corollary}{0}
\setcounter{Lemma}{0}
\setcounter{Definition}{0}
\setcounter{Remark}{0}

\subsection{Basic properties}
As a first step, we need to establish the  following easy properties of
$P^n_k$ as a function of the global constraints
$Q=(Q_{\min},Q_{\max})$.

\ms
\ni {\bf P1.} {\em For every $k\!\in \{0,\ldots,n\}$, and every couple of global constraints (at time $k$),
\[
P^n_k(Q,(V_k,\F_k)_{0\le k\le n-1})= P_0^{n\!-\!k}(Q,(V_{k+\ell},\F_{k+\ell})_{0\le \ell\le
n\!-\!k-1}).
\]
}
This is obvious from~(\ref{(2.3)}).

\ms
\ni {\bf P2.} {\em If $Q$ and $Q'$ are two admissible global constraints (at time $k$) then
\[
P^n_k(Q)=P^n_k(Q') \quad \mbox{ on the event }\quad \{Q=Q'\}.
\]
}
\ni {\bf Proof.} Owing to {\bf P1}  one may assume  without loss of generality that $k=0$. Let $q$ and $ q'$ be two admissible
controls with respect to $Q$ and $ Q'$ respectively. Set $\tilde q_\ell = q_\ell \mbox{\bf 1}_{\{Q=  Q'\}} +   q'_\ell\mbox{\bf
1}_{\{Q\neq  Q'\}}$, $\ell=0,\ldots,n-1$. The control  $\tilde q$ is admissible with respect to $ Q'$ since
$\{Q=  Q'\}\!\in {\cal F}_0$. Furthermore,
\begin{eqnarray*}
\mbox{\bf 1}_{\{Q=\  Q'\}} \E\left(\sum_{\ell=0}^{n-1} q_\ell V_\ell\,|\, {\cal F}_0\right)\!&\!=\!&\!\E\left(\sum_{\ell=0}^{n-1}
\mbox{\bf 1}_{\{Q= Q'\}} q_\ell V_\ell\,|\, {\cal F}_0\right)=  \mbox{\bf 1}_{\{Q= Q'\}}
\E\left(\sum_{\ell=0}^{n-1}  \widetilde q_\ell V_\ell\,|\, {\cal F}_0\right).\\
&\le&  \mbox{\bf 1}_{\{Q=  Q'\}} P^n_0(  Q') \quad a.s.
\end{eqnarray*}
Hence,
\[
\mbox{\bf 1}_{\{Q= Q'\}}P^n_0(Q) = \supess \left\{\mbox{\bf 1}_{\{Q= Q'\}}
\E\left(\sum_{\ell=0}^{n-1} q_\ell V_\ell\,|\, {\cal F}_0\right),\;
q\hbox{
$Q$-admissible} \right \} \le  \mbox{\bf 1}_{\{Q= Q'\}} P^n_0( Q') \quad
a.s.
\]
The equality follows by symmetry.$\cqfd$

\ms
\ni {\bf P3.}   {\em Let $k\!\in \{0,\ldots,n-1\}$. The  set of admissible global constraints $Q= (Q_{\min},Q_{\max})$ at time $k$
is convex and the mapping
$Q\mapsto P^n_k(Q)$ is concave in the following sense: if $Q$ and
$Q'$ are two couples of admissible constraints, then for  every random variable $\lambda:(\Omega,{\cal F}_k)\to [0,1]$, $\lambda
Q+(1-\lambda)Q'$ is an admissible couple of constraints and
\[
P^n_k(\lambda
Q+(1-\lambda)Q') \ge \lambda P^n_k(Q)+(1-\lambda)P^n_k(Q')\quad a.s.
\]
 Furthermore $ Q_{\min}\mapsto P^n_k(Q_{\min},Q_{\max})$ is
non-increasing and $ Q_{\max}\mapsto P^n_k(Q_{\min},Q_{\max})$ is
non-decreasing   $i.e.$
\[
P^n_k(Q)\le P^n_k(Q') \qquad  a.s. \qquad \mbox{\rm  on }
\quad  \{Q'_{\min}\le Q_{\min}\le Q_{\max}\le Q'_{\max}\}.
\]
}
\ni{\bf Proof.} One may assume  by {\bf P1}  that $k=0$. The convexity of admissible global constraints is obvious.
As concerns the concavity of the value function, note that if $q$ and $q'$ are locally admissible controls then $\lambda
q+(1-\lambda) q':=(\lambda q_k+(1-\lambda)q'_k)_{0\le k\le n-1}$ is still locally admissible. If $q$ and $q'$ satisfy the  global
constraints induced by $Q$ and
$Q'$ respectively, then $\lambda q+(1-\lambda) q'$ always satisfies that induced by  $\lambda
Q+(1-\lambda)Q'$. Consequently, using that $\lambda$ is ${\cal F}_0$-measurable,
\begin{eqnarray*}
&&P^n_0(\lambda
Q+(1-\lambda)Q')\\
& \ge&\supess\left\{\E(\sum_{k=0}^{n-1}(\lambda
q_k+(1-\lambda)q'_k)V_k\,|\,{\cal F}_0),\, q,\, q' \hbox{ locally admissible}, \, \bar q_{n}\!\in
[Q_{\min},Q_{\max}],\,\bar q'_{n}\!\in [Q'_{\min},Q'_{\max}] \right\}\\
\!&\!=\!&\! \lambda\, \supess\left\{\E(\sum_{k=0}^{n-1}
q_kV_k\,|\,{\cal F}_0), \, q \hbox{ locally admissible}, \, \bar q_{n}\!\in [Q_{\min},Q_{\max}] \right\}\\
&&+(1-\lambda)\,\supess\left\{\E(\sum_{k=0}^{n-1}q'_kV_k\,|\,{\cal
F}_0), \, q' \hbox{ locally admissible}, \,\bar q'_{n}\!\in
[Q'_{\min},Q'_{\max}] \right\}\\
\!&\!=\!&\! \lambda P^n_0(Q)+ (1-\lambda)P^n_0(Q').
\end{eqnarray*}

The monotony property is as follows: let $A:= \{Q'_{\min}\le Q_{\min}\le
Q_{\max}\le Q'_{\max}\}\!\in \F_0$ and let $\tilde q'$ be a fixed
$Q'$-admissible control. Then, for every $Q$-admissible  control $q$, set
\[
q':= q\mbox{\bf 1}_{_A}+ \tilde q'\mbox{\bf 1}_{^cA}.
\]
Then $q'$ is clearly $Q'$-admissible and
\begin{eqnarray*}
\mbox{\bf 1}_{_A}\E\left(\sum_{k=0}^{n-1}q_kV_k\,|\,{\cal F}_0\right)&=&\E\left(\sum_{k=0}^{n-1}\mbox{\bf 1}_{_A}q_kV_k\,|\,{\cal
F}_0\right)\\&=&
\E\left(\mbox{\bf 1}_{_A} \sum_{k=0}^{n-1}q'_kV_k\,|\,{\cal F}_0\right)\\
&\le& \mbox{\bf
1}_{_A}P_0^n(Q')\qquad a.s.\\
 \mbox{ so that }\hskip 5 cm \mbox{\bf 1}_{_A}P^n_0(Q)&\le&
\mbox{\bf 1}_{_A}P^n_0(Q')\qquad a.s.\hskip 2,25 cm \qquad\cqfd
\end{eqnarray*}

\bs
\ni {\bf P4.} {\em Let $k\!\in \{0,\ldots,n-1\}$. Let  $Q^{(n)}$, $n\ge
1$, be a sequence of admissible global constraints such that
$Q^{(n)}_{\max}\to 0$ as $n\to \infty$. Then
\[
P_k^n(Q^{(n)})\longrightarrow 0\qquad \mbox{ as }\qquad n\to \infty .
\]
}
\ss
\ni{\bf Proof.} Owing to {\bf P1}, one may assume $k=0$ without loss of
generality. The result is straightforward once noticed that for every
$Q$-admissible control
\[
\left|\E\left(\sum_{k=0}^{n-1} q_k V_k\,|\, \F_0\right)\right|  \le
Q_{\max}
\E\left(\max_{0\le k\le n-1} |V_k|\,|\, \F_0\right)  .\qquad\cqfd
\]

\bs
\ni {\bf P5.} {\em Let $T^+(n):=\{(u,v)\!\in\R_+^2,\, 0\le u\le v\le n\}$ and let $k\!\in \{0,\ldots,n-1\}$. There is   process
$\Pi^n_k:(T^+(n)\times \Omega, {\cal B}or(T^+(n))\otimes {\cal F}_k)\to \R$
such that

\ss
\ni$(i)$ $(u,v)\mapsto \Pi^n_k(u,v,\omega)$ is concave and continuous on $T^+(n)$ for every $\omega\!\in \Omega$,

\ss \ni $(ii)$ $u\mapsto \Pi^n_k(u,v,\omega)$ is non-increasing on
$[0,v]$ and  $v\mapsto \Pi^n_k(u,v,\omega)$ is nondecreasing on
$[u,n]$ for every $\omega\!\in \Omega$.

\ss
\ni $(iii)$ For every admissible constraint $Q=(Q_{\min},Q_{\max})$ at  time $k$,
$P^n_k(Q)(\omega)= \Pi^n_k(Q(\omega), \omega)$
$\P(d\omega)$-$a.s.$}

\ss
\ni{\bf Proof.} Classical consequence of {\bf P3}:  for every
$(r,s)\!\in T^+(n)\cap \Q^2$ set
\[
\Pi^n_k((r,s),\omega):= P^n_k((r,s))(\omega).
\]
Then, set for every $(u,v)\!\in T^+(n)$
\[
\Pi^n_k((u,v),\omega):=a.s.\lim_{(r,s)\to (u,v),  (r,s)\in T^+(n)\cap \Q^2,
r\le u\le v\le s} P^n_k((r,s))(\omega).
\]

One shows using the concavity and monotony properties established in {\bf
P3} that the above limit does exist and that , $\P(d\omega)$-$a.s.$
$(u,v)\mapsto
\Pi^n_k((u,v),\omega)$ is continuous on $T^+(n)$  and that
$P^n_k(Q)(\omega)=\Pi^n_k( Q(\omega),\omega)$
$\P$-$a.s.$.


\bs
 \ni {\bf P6.}  Let $Q=(Q_{\min}, Q_{\max})$ be a couple of global
constraints (at time $0$).
\[
P^n_0(Q)= \supess\left\{\E\left(\sum_{\ell=0}^{n-1} q_\ell(V_\ell)^+\,|\,{\cal F}_0\right) ,\; q \hbox{
$({\F},Q)$-admissible}\right\} \;\mbox{ on the event }\;
\{Q_{\min}=0\} \!\in {\cal F}_0.
\]

\ni{\bf Proof.} This   follows from the simple remark that one can define from any $({\F},Q)$-admissible control $q$ a
new $({\F},Q)$-admissible control $\widetilde q$ by $\widetilde q_\ell := q_\ell\mbox{\bf 1}_{\{V_\ell\ge 0\}\cap \{Q_{\min}=0\}
}+q_\ell\mbox{\bf 1}_{\{Q_{\min}>0\}
}$ and
\[
\mbox{\bf 1}_{\{Q_{\min}=0\}}\sum_{k=0}^{n-1} q_\ell V_\ell \le \mbox{\bf
1}_{\{Q_{\min}=0\}}\sum_{k=0}^{n-1} q_\ell(V_\ell)^+= \mbox{\bf
1}_{\{Q_{\min}=0\}}\sum_{k=0}^{n-1}
\tilde q_\ell (V_\ell)^+. \cqfd
\]

\subsection{Dynamic programming principle}
 The main consequence of  {\bf P5} is  that   at every time $k$, one may assume without loss of generality that  the couple of
admissible global constraints
$Q=(Q_{\min},Q_{\max})$ is deterministic since, for any possibly random admissible constraint $Q$ (at time $k$) and every $\omega\!\in
\Omega$, $
 P^n_k(Q)(\omega)= \Pi^n_k(\omega, Q(\omega))$.

As a consequence, for notational convenience, we will still denote $P^n_k$ instead of $\Pi^n_k$ so that for any admissible global
constraints $Q$ at time $k$
\[
P^n_k(Q)(\omega)= P^n_k(\omega,x)_{|x=Q(\omega)}.
\]

\begin{Thm}(Backward Dynamic Programming Principle) Set $P^n_{n}\equiv 0$.

\noindent $(a)$ {\sc Local Dynamic programming formula.} For every $k\!\in \{0,\ldots,n-1\}$ and every couple
$Q=(Q_{\min},Q_{\max})$ of deterministic admissible global constraints at time
$k$
\begin{equation}\label{PgmDyn}
\hskip -0.35 cm P^n_k(Q) = \sup \left\{xV_k + \E\left(P^n_{k+1}(\chi^{n\!-\!k-1}(Q,x))\,|\,{\cal F}_k\right),\, x\!\in
I^{n-1-k}_Q\right\}
\end{equation}
where $\chi^{M}(Q,x)= ((Q_{\min}-x)^+,
(Q_{\max}-x)\!\wedge\!M)$ and $I^M_Q:=[(Q_{\min}-M)^+\!\wedge\! 1,Q_{\max}\!\wedge\! 1]$.

\medskip
\noindent  $(b)$  {\sc Global Dynamic programming formula.} For every couple
$Q=(Q_{\min},Q_{\max})$ of   admissible global constraints at time $0$,
the price of the contract at time $k\!\in\{0,\ldots,n-1\}$ is given by   $P^n_k(Q^{k,*})$  where
\begin{eqnarray}\label{PgmDynGlob}
  Q^{k,*} &:= &(Q_{\max}-\bar q^*_k,Q_{\min}-\bar q^*_k)\hskip 1 cm \mbox{(residual global constraints) } \\
\nonumber \mbox{with } \hskip 0.8 cm q^*_k &=&  q^*_k(Q^{k,*}) \hskip 5 cm\\
 \label{Optiqk}
 q^*_k(Q)&:=&{\rm argmax}_{x\in I^{n-1-k}_Q}\left(xV_k + \E\left(P^n_{k+1}(\chi^{n-k-1}(Q,x))\,|\,{\cal F}_k\right)\right),\;
k=0,\ldots,n-1.\qquad
\end{eqnarray}
Furthermore,
\[
P^n_k(Q^{k,*}) = \E\left(\sum_{\ell=k}^{n-1} q^*_\ell V_\ell\,|\, {\cal F}_\ell\right).
\]
\end{Thm}

\noindent {\bf Remark.} The definition (\ref{Optiqk}) may be
ambiguous when argmax is not reduced to a single point. Then, one
considers $\min\,$argmax to define $q^*_k(Q)$.

\bigskip
\ni {\bf Proof.} $(a)$ It is clear that owing to {\bf P1}  the case $k=0$ is the only one to be proved. As a first step, we prove
that
\begin{equation}\label{PgmDyn2}
\hskip -0.35 cm P^n_0(Q) = \supess\left\{q_0V_0 + \E\left(P^n_{1}(\chi^{n-1}(Q,q_0))\,|\,{\cal F}_0\right),\,
q_0:(\Omega,{\cal F}_0)\to I^{n-1}_Q\right\}
\end{equation}

\ms
\ni \fbox{$\le$} : Let $q=(q_k)_{0\le k\le n-1}$ be a $Q$-admissible control. Then, $q_0$ is $[0,1]$-valued (as well as the $q_k$'s) and
\[
  Q_{\min}-(n-1)\le Q_{\min}-(q_1+\cdots+q_{n-1})  \le q_0\le     Q_{\max}-(q_1+\cdots+q_{n-1})\le Q_{\max}
\]
so that  $q_0$ is $I^{n-1}_Q$-valued. Furthermore,
\begin{equation}\label{IneqRecContr}
\max(0, Q_{\min}-q_0) \le q_1+\cdots+q_{n-1}\le   (Q_{\max}-q_0)\!\wedge\!(n-1).
\end{equation}
 note that
$\chi^{n-1}(Q,q_0)$ is an admissible couple of ($\F_0$-measurable) constraints at time $1$. Consequently,
\begin{eqnarray*}
\E\left(\sum_{1\le \ell\le n-1} q_\ell V_\ell\,|\,{\cal F}_0\right)\!&\!=\!&\!\E\left(\E\left(\sum_{1\le \ell\le n-1} q_\ell V_\ell\,|\,{\cal
F}_{1}\right)\,|\,{\cal F}_0\right)\\
&\le& \E\left(P^n_1(\chi^{n-1}(Q,q_0))\,|\,{\cal F}_0\right)\quad a.s.
\end{eqnarray*}
where the last inequality follows from the definition of $P^n_1$, $\chi^{n-1}$, (\ref{IneqRecContr}) and the monotony of conditional
expectation.  Then,
\[
\E\left( \sum_{0\le \ell\le n-1} q_\ell V_\ell\,|\,{\cal F}_0\right)\le
\E\left(q_0V_0+P^n_1(\chi^{n-1}(Q,q_0))\,|\,{\cal F}_0\right)\quad a.s.
\]
One concludes that
\[
P^n_0(Q) \le \supess\left\{\E\left(q_0V_0+P^n_1(\chi^{n-1}(Q,q_0))\,|\,{\cal
F}_0\right),\; q_0:(\Omega,{\cal F}_0)\to I^{n-1}_Q\right\}.
\]

\ms
\ni \fbox{$\ge$} : We  proceed as usual by proving a bifurcation property
for the controls.  Let $q_0$ and $ q'_0$ be two  $I^{n-1}_Q$-valued
$\F_0$-measurable random variables. Set
\[
A_0:= \left\{ q_0 V_0 +
\E\left(P^n_{1}(\chi^{n-1}(Q,q_0))\,|\,{\cal F}_0\right) >   q'_0
V_0 + \E\left(P^n_{1}(\chi^{n-1}(Q,  q'_0))\,|\,{\cal F}_0\right)
\right\}\!\in {\cal F}_0.
\]
and
\[
\tilde q_0= q_0\mbox{\bf 1}_{A_0} +  q'_0 \mbox{\bf 1}_{^c A_0}:
(\Omega,{\cal F}_0)\to I^{n-1}_Q .
\]
Then,  one checks that
\[
\tilde q_0 V_0 + \E\left(P^n_{1}(\chi^{n-1}(Q,\tilde q_0))|{\cal
F}_0\right)=
\max_{y=q_0, q'_0}  \left(yV_0 +
\E\left(P^n_{1}(\chi^{n-1}(Q,y))|{\cal F}_0\right)\right).
\]

Consequently there exists  a sequence  $q_0^{(n)}$ of $[0,Q_{\max}\wedge 1]$-valued random variables such that
\begin{eqnarray*}
\supess\left\{ q_0 V_0 + \E\left(P^n_{1}(\chi^{n-1}(Q,  q_0))\,|\,{\cal F}_0\right),\,
q_0:(\Omega,{\cal F}_0)\to I^{n-1}_Q\right\}&&\\
= \sup_n q^{(n)}_0V_0 + \E\left( P^n_{1}(\chi^{n-1}(Q,  q^{(n)}_0))\,|\,{\cal
F}_0\right) &&\\
 = \lim^{\;\qquad_\uparrow}_n q^{(n)}_0V_0 + \E\left(
P^n_{1}(\chi^{n-1}(Q,  q^{(n)}_0))\,|\,{\cal F}_0\right).
\end{eqnarray*}
One may assume by applying the above bifurcation property that the above supremum holds as a nondecreasing limit as $n\to \infty$.

Now, for every fixed $n\ge 1$, there exists a sequence of $[0,1]^{n-1}$-valued random vectors $(q_k^{(n,m)})_{1\le k\le n-1}$ such
that $\sum_{k=1}^{n-1} q_k^{(n,m)}\!\in [(Q_{\min}-q^{(n)}_0)^+, Q_{\max}-q^{(n)}_0]$ and
\[
 P^n_{1}(\chi^{n-1}(Q,  q^{(n)}_0)) = \lim^{\;\quad_\uparrow}_m \left(\sum_{k=1}^{n-1}
q_k^{(n,m)} V_k\,|\, {\cal F}_1\right)\quad a.s.
\]
 where we used that the admissible  sequences $(q_k)_{1\le k\le n-1}$ for the problem starting at $1$ clearly satisfy the bifurcation principle
due to the homogeneity of conditional expectation $\E(\,.\,|\,{\cal F}_0)$ with respect to ${\cal F}_0$ -measurable r.v.. Consequently,
\begin{eqnarray*}
 q^{(n)}_0V_0+P^n_{1}(\chi^{n-1}(Q,  q^{(n)}_0)) \!&\!=\!&\!
\lim^{\;\qquad_\uparrow}_m q^{(n)}_0V_0+\E\left(\sum_{k=1}^{n-1}
q_k^{(n,m)} V_k\,|\, {\cal F}_1\right)\quad a.s.\\
\E\left( q^{(n)}_0V_0+P^n_{1}(\chi^{n-1}(Q,  q^{(n)}_0))\,|\, {\cal
F}_0\right)\!&\!=\!&\!\lim^{\;\qquad_\uparrow}_m\E\left(
q^{(n)}_0V_0+\sum_{k=1}^{n-1} q_k^{(n,m)} V_k\,|\, {\cal F}_0\right)\quad
a.s.
\end{eqnarray*}
where we used the conditional Beppo Levi Theorem. Note that for every $n,\, m\ge 1$, $(q^{(n)}_0, q_k^{(n,m)}, k=1,\ldots,n-1)$ is an admissible
control with respect to $Q$ since $x+(Q_{\min}-x)^+\ge Q_{\min}$. Hence
\begin{eqnarray*}
\E\left( q^{(n)}_0V_0+P^n_{1}(\chi^{n-1}(Q,  q^{(n)}_0))\,|\, {\cal F}_0\right)&\le
&\supess\left\{\!\E\left( q_0V_0\!+\!\sum_{k=1}^{n\!-\!1} q_k V_k\,|\, {\cal F}_0\right), \;q \, \hbox{ $Q$-admissible}\right\}\;
a.s.\\
\!&\!=\!&\! P^n_0(Q) \quad a.s.
\end{eqnarray*}

To pass from~(\ref{PgmDyn2}) to~(\ref{PgmDyn}) is standard using {\bf P5}. Let $q_0:(\Omega,{\cal F}_0)\to
I^{n-1}_Q$.
\begin{eqnarray*}
q_0V_0+\E(P^n_1(\chi^{n-1}(Q,q_0))\,|\,\F_0)&=&
q_0V_0+\E(P^n_1(y)\,|\,\F_0)_{|y=\chi^{n-1}(Q,q_0)}\\&
\le& \sup_{x\in I^{n-1}_Q}(xV_0+\E(P^n_1(y)\,|\,\F_0)_{|y=\chi^{n-1}(Q,x)})\\
&=& \sup_{x\in I^{n-1}_Q}(xV_0+\E(P^n_1(\chi^{n-1}(Q,x))\,|\,\F_0).
\end{eqnarray*}
Conversely, setting $q^x_0(\omega):=x\!\in I^{n-1}_Q$  yields the reverse inequality.

\medskip
\noindent $(b)$ This item follows from {\bf P5}. $\quad\cqfd$

\section{Affine value function with bang-bang optimal controls}\label{piecaffine}

\setcounter{equation}{0}
\setcounter{Assumption}{0}
\setcounter{Theorem}{0}
\setcounter{Proposition}{0}
\setcounter{Corollary}{0}
\setcounter{Lemma}{0}
\setcounter{Definition}{0}
\setcounter{Remark}{0}

\subsection{The main result}
We recall that, for every integer $n\ge 1$, the triangular set of
admissible values for a couple of (deterministic) global constraint
(at time $0$), as introduced in {\bf P5}, is defined as:
$$
T^+(n):=\{(u,v),\; 0\le u\le v\le n\}.
$$
Then we will define a triangular tiling of $T^+(n)$ as follows: for every couple of integers $(i,j)$, $0\le i\le j\le n-1$,
\[
T^+_{i\,j}:= \{(u,v)\!\in [i,i+1]\!\times\![j,j+1],\; v\ge u+j-i\}\quad \mbox{and}\quad T^-_{i\,j}:= \{(u,v)\!\in
[i,i+1]\!\times\![j,j+1],\; v\le u+j-i\}.
\]
One checks that
\[
T^+(n) = (\bigcup_{0\le i\le j\le n-1} T^+_{i\,j}) \bigcup (\bigcup_{0\le i< j\le n-1}T^-_{i\,j}).
\]

 \begin{Thm}\label{main} The multi-period swing option premium with deterministic global constraints $Q:=(Q_{\min}, Q_{\max})\!\in
T^+(n)$  as defined by~(\ref{(2.3)})  is always obtained as the result of an optimal strategy.

\ms
\ni $(a)$ The value function (premium):

\ss -- the mapping $Q\mapsto P_0(Q,\F)$ is  a concave, continuous,
piecewise affine process, affine  on every triangle $T^\pm_{i,j}$
of the tiling of $T^+(n)$. Furthermore,
\[
P^n_0(0,0)=0,\qquad P^n_0(0,n)=
\E\left(\sum_{k=0}^{n-1} V_k^+\,|\,\F_0\right)\quad \mbox{ and } \quad
 P^n_0(n,n)=
\E\left(\sum_{k=0}^{n-1}V_k\,|\,\F_0\right).
\]

\ss -- If $V_i \ge 0$ $a.s.$ for every $i\!\in\{0,\ldots,n-1\}$,
then,  for every $Q=(Q_{\min}, Q_{\max})\!\in T^+(n)$,
$P_0(Q,\F)=P_0((0, Q_{\max}),\F)= P_0((Q_{\max}, Q_{\max}),\F)$
$a.s.$. (in particular $Q_{\max}\mapsto P_0((Q_{\max},
Q_{\max}),\F)$ is $a.s.$ non-decreasing).

\ms
\ni $(b)$ The optimal control:

\ss
-- If the  global constraint $Q=(Q_{\min},Q_{\max}) \!\in \N^2\cap T^+(n)$,  then   there always
exists a bang-bang optimal control $q^*=(q^*_k)_{0\le k\le n-1}$ with $q^*_k$ is $\{0,1\}$-valued for every $k=0,\ldots,n-1$. .

\ss
-- If $V_i \ge 0$ $a.s.$ for every $i\!\in\{0,\ldots,n-1\}$, then      there
always exists a bang-bang optimal control which satisfies
$\sum_{0\le k\le n-1} q^*_k =Q_{\max}$.

\ss
-- Otherwise the optimal control is not bang-bang as emphasized by the case $n=2$ (see proposition~\ref{n=2} below)
 \end{Thm}

We will first inspect the case of a two period swing contract. It will illustrate in a simpler
setting the approach developed in the general case. Furthermore, we will obtain a slightly more
precise result about the optimal controls.

\subsection{The two period option} We assume $n=2$ throughout this
section.  The  first result is the following
\begin{Pro} Let $Q=(Q_{\min},Q_{\max})\!\in T^+(2)$ denote an admissible global
constraint and $I^{1}_Q:=[(Q_{\min}-1)^+, Q_{\max}\!\wedge\!1]$. There is an optimal control
$q^*=(q^*_0,q^*_1)$  given by
\begin{eqnarray}
\label{q*1} q^*_0\!&\!=\!&\!\displaystyle{{\rm argmax}_{x\in I^1_Q}} \left\{x
V_0+1\!\wedge\! (Q_{\max}\!-x )\E (V_1^+|\F_0)-
(Q_{\min}\!-x )^+\E(V^-_1|\F_0)\right\}\\
\label{q*2}q^*_1\!&\!=\!&\! 1\!\wedge\! (Q_{\max}\!-q^*_0)\mbox{\bf 1}_{\{V_1\ge 0\}}+
(Q_{\min}\!-q^*_0)^+\mbox{\bf 1}_{\{V_1< 0\}}
\end{eqnarray}
so that
\[
P^2_0(Q,\F) = \E ( q^*_0V_0+q^*_1V_1\,|\, \F_0).
\]
\end{Pro}

\noindent {\bf Proof.} Let $q=(q_0,q_1)$ be an admissible control:
$q_0+q_1\!\in [Q_{\min},Q_{\max}]$ and $q_i$ are $[0,1]$-valued $\F_i$-measurable, $i=0,1$.
Consequently $q_0$ is $I^1_Q$-valued and
$q_1$ is $[(Q_{\min}-q_0)^+, (Q_{\max}-q_0)\!\wedge\!1]$-valued. Hence
\begin{equation}\label{ineq1}
q_0V_0+q_1V_1 \le q_0 V_0 + 1\!\wedge\! (Q_{\max}\!-q_0) V_1^+ -
(Q_{\min}\!-q_0)^+ V^-_1.
\end{equation}
On the other hand
\begin{eqnarray*}
\E\left(q_0 V_0 + 1\!\wedge\! (Q_{\max}\!-q_0) V_1^+\right.  \!\! \!\!&\!\! -\!\! & \!\!\!\!\left.
(Q_{\min}\!-q_0)^+ V^-_1|\F_0\right)\\
\!&\!=\!&\! q_0
V_0+1\!\wedge\! (Q_{\max}\!-q_0)\E (V_1^+|\F_0)-
(Q_{\min}\!-q_0)^+\E(V^-_1|\F_0).
\end{eqnarray*}
The mapping $x\mapsto x\, V_0+1\!\wedge\! (Q_{\max}\!-x)\E
(V_1^+|\F_0)- (Q_{\min}\!-x)^+\E(V^-_1|\F_0) $ (called the {\em
objective variable} from now on) is piecewise affine on $I^1_Q$
with $\F_0$-measurable coefficients so the above definition of
$q^*_0$ defines an $\F_0$-measurable $I^1_Q$-valued random
variable. Now, combining the above inequalities yields
\begin{eqnarray*}
\E(q_0V_0+q_1V_1 \,|\,\F_0)&\le & q_0V_0+1\!\wedge\!(Q_{\max}-q_0) \E(V_1^+\,|\,\F_0)- (Q_{\min}\!-q_0)^+
\E(V^-_1|\F_0) \\
&\le &  \sup_{q_0\in I^1_Q} q_0V_0+1\!\wedge\!(Q_{\max}-q_0) \E(V_1^+\,|\,\F_0)- (Q_{\min}\!-q_0)^+
\E(V^-_1|\F_0)\\
&=& \E(q^*_0V_0+q^*_1V_1\,|\,\F_0).\quad
\cqfd
\end{eqnarray*}

In the proposition below we investigate in full details the case $n=2$.

\begin{Pro}\label{n=2} Let $n=2$. The two period swing option premium with admissible global constraints
$Q=(Q_{\min}, Q_{\max})\!\in T^+(2)$ as defined by~(\ref{swing})  is always obtained as the result of an optimal strategy.

\ms
\ni $(a)$ The optimal control:

\ss
-- If the global constraints $Q_{\min}, \,Q_{\max}$ only take integral values (in $\{0,1,2\}$)
then there always exists a
$\{0,1\}$-valued bang-bang optimal control. When $Q$ simply satisfies  $Q_{\max}-Q_{\min}\!\in\{0,1,2\}$,
there always exists a bang-bang
optimal control.

\ss
-- If $V_0,\, V_1\ge 0$ $a.s.$, then   any  optimal
 control  $q^*$ is $a.s.$ bang-bang and satisfies $ q^*_0+q^*_1= Q_{\max}$ on $\{V_i>0,\,
i=1,2\}$.

\ss -- Otherwise the optimal control is generally not bang-bang.

\ms
\ni $(b)$ The value function (premium):

\ss
-- the mapping $Q=(Q_{\min},Q_{\max})\mapsto P_0(Q, {\cal F})$ is  affine on the
four triangles $T^\pm_{i,j}$ that tile $T^+(2)$.

\ss
-- Furthermore, when $V_0$ and $V_1$ are $a.s.$ non negative,
\[
P^2_0(Q,\F) =  (Q_{\max}\!-\!1)^+\!\wedge\! 1 \left(
V_0\!\wedge\!\E(V_1|\F_0)\right)+ (Q_{\max}\!\wedge\!1)(V_0\vee
\E(V_1|\F_0)).
\]
\end{Pro}

The objective variable being  piecewise affine on $I^1_Q$,    $q^*_0$  is   equal either to  one of its monotony
breaks or to the
endpoints of $I^1_Q$. Consequently, a careful inspection of all  possible situations for the global constraints
yields the complete set
of    explicit optimal rules for the optimal exercise of the swing option involving the values $V_0$ and
$\E(V_1^\pm|\F_0)$ (expected gain
or loss  at time
$0$) at time $0$ and $V_1$ at time $1$.

\bs
\noindent \fbox{$Q\!\in T^+_{00}$ $i.e.$ $Q_{\min} \le
Q_{\max}\le 1$: } $I^1_Q=[0,Q_{\max}]$  and
  the objective variable reads
\[
q_0V_0+ (Q_{\max}\!-q_0) \E(V_1^+|\F_0) - (Q_{\min}\!-q_0)^+\E(V_1^-|\F_0)
\]
with  one monotony break at $Q_{\min}$. One checks that
\[
\begin{array}{lll}
q^*_0=Q_{\max},& q^*_1=0 &\mbox{on }\; \{V_0\ge \E(V_1^+|\F_0) \},\\
q^*_0= Q_{\min},& q^*_1= (Q_{\max}\!-Q_{\min})\mbox{\bf 1}_{\{V_1\ge 0\}} &\mbox{on }\;\{\E(V_1|\F_0)
\le V_0<\E(V_1^+|\F_0)
\},\\
 q^*_0= 0, &  q^*_1 = Q_{\max} \mbox{\bf 1}_{\{V_1\ge 0\}}+ Q_{\min} \mbox{\bf 1}_{\{V_1<0\}}
&\mbox{on }\{V_0 < \E(V_1|\F_0) \}.
\end{array}
\]

Note that
\[
q^*_0= Q_{\min},\;  q^*_1= Q_{\max}\!-Q_{\min} \;\mbox{ on }\;\{\E(V_1|\F_0)
\le V_0<\E(V_1^+|\F_0)\}\cap \{V_1\ge 0\}
\]
so that the above optimal control is not bang-bang on this event except if $Q_{\min}\!\in\{0,1\}$ or $Q_{\max}=Q_{\min}$.

\bs
\noindent \fbox{$Q\!\in T_{01}^+$ $i.e.$ $Q_{\min} \le 1\le Q_{\max}\le 2$, $Q_{\min} \le
Q_{\max}\!-1$: }  $I^1_Q= [0,1]$ and the objective variable reads
\[
q_0V_0+1\!\wedge\!(Q_{\max}\!-q_0) \E(V_1^+|\F_0) - (Q_{\min}\!-q_0)^+\E(V_1^-|\F_0)
\]
with  monotony breaks
 at $Q_{\min}$, $Q_{\max}\!-1$. One checks that
\[
\begin{array}{lll}
q^*_0=1,& q^*_1= (Q_{\max}\!-1)\mbox{\bf 1}_{\{V_1\ge 0\}} &\mbox{on }\; \{V_0\ge \E(V_1^+|\F_0) \},\\
q^*_0= Q_{\max}\!-1,& q^*_1= \mbox{\bf 1}_{\{V_1\ge 0\}} &\mbox{on
}\;\{0 \le V_0<\E(V_1^+|\F_0)
\},\\
 q^*_0= Q_{\min}, &  q^*_1 =  \mbox{\bf 1}_{\{V_1\ge 0\}}   &\mbox{on }\{ -\E(V_1^-|\F_0)  \le
V_0<0\}\\   q^*_0= 0, &  q^*_1 =   \mbox{\bf 1}_{\{V_1\ge 0\}}+ Q_{\min} \mbox{\bf
1}_{\{V_1< 0\}}  &\mbox{on }\{ V_0<-\E(V_1^-|\F_0) \}.
\end{array}
\]

Note that
$$
q^*_0+q^*_1 = Q_{\max}\!-1\quad \mbox{ on } \quad\{0 \le V_0<\E(V_1^+|\F_0)\}\cap \{V_1<0\}
$$
so that the control is not bang-bang on this event, except if  $Q_{\max}\in\{1,2\}$  or  $Q_{\max}=  1+Q_{\min}$, since the local
control
$q^*_0$ and the global constraint are not saturated.  Likewise
$$
q^*_0+q^*_1 = 1+Q_{\min}\quad \mbox{ on } \quad\{ -\E(V_1^-|\F_0)  \le
V_0<0\}\cap \{V_1>0\}
$$
and the optimal control is not bang-bang on this event, except when $Q_{\min}\!\in\{0,1\}$ or  $Q_{\max}=  1+Q_{\min}$.

Note that {\em both events correspond to prediction errors}: $V_1$ has not the predicted sign. Moreover, these events are $a.s.$  empty
when $V_i\ge 0$ $a.s.$, $i=1,2$. On all  other events the  optimal control is bang-bang.

\bs
\noindent \fbox{$Q\!\in T_{01}^-$ $i.e.$ $Q_{\min} \le 1\le Q_{\max}\le2$, $Q_{\min} \ge Q_{\max}\!-1$:}
  Then the monotony breaks of the objective process (with the same expression as in the former case) still
are $Q_{\min}$, $Q_{\max}\!-1$. A careful inspection of the four possible cases leads to
\[
\begin{array}{lll}
q^*_0=1,& q^*_1= (Q_{\max}\!-1)\mbox{\bf 1}_{\{V_1\ge 0\}} &\mbox{on }\; \{V_0\ge \E(V_1^+|\F_0) \},\\
q^*_0= Q_{\min},& q^*_1= (Q_{\max}\!-Q_{\min})\mbox{\bf 1}_{\{V_1\ge 0\}} &\mbox{on
}\;\{\E(V_1|\F_0) \le V_0<\E(V_1^+|\F_0)
\},\\
 q^*_0= Q_{\max}\!-1, &  q^*_1 =  \mbox{\bf 1}_{\{V_1\ge 0\}} + (Q_{\min}\!-Q_{\max}+1) \mbox{\bf
1}_{\{V_1< 0\}}   &\mbox{on }\{ -\E(V_1^-|\F_0)  \le
V_0<\E(V_1|\F_0)\}\\   q^*_0= 0, &  q^*_1 =  \mbox{\bf 1}_{\{V_1\ge 0\}}+ Q_{\min} \mbox{\bf
1}_{\{V_1< 0\}}  &\mbox{on }\{ V_0<-\E(V_1^-|\F_0) \}.
\end{array}
\]
Note that  on the event
$$
\{ -\E(V_1^-|\F_0)  \le
V_0<\E(V_1|\F_0)\}\cap \{V_1<0\}
$$
the optimal control is not bang-bang, except if
$Q_{\max}\!\in\{1,2\}$ or $Q_{\max}=Q_{\min}$ (both  $q^*_0$ and $q^*_1$ are
$(0,1)$-valued) or $Q_{\max}= Q_{\min}+1$ ($q^*_0=Q_{\min}$, $q^*_1=0$);   and on the event
$$
\{\E(V_1|\F_0) \le V_0<\E(V_1^+|\F_0)
\}\cap  \{V_1>0\}
$$
the optimal control is not bang-bang  either (except if
$Q_{\min}\!\in\{0,1\}$ or $Q_{\max}=Q_{\min}$ or $Q_{\max}= Q_{\min}+1$) by similar arguments.

Note that  these  events {\em do not correspond to an error of prediction}. On all other events
the  optimal control is bang-bang.

\bs
\noindent \fbox{$Q\!\in T_{11}^+$ $i.e.$ $1<Q_{\min}  \le Q_{\max}\le 2$:} The objective variable is defined on $I^1_Q=[
Q_{\min}\!-1,1]$ by
\[
q_0V_0+1\!\wedge\!(Q_{\max}\!-q_0) \E(V_1^+|\F_0) -(Q_{\min}\!-q_0)\E(V_1^-|\F_0)
\]
with  only   one breakpoint at   $Q_{\max}\!-1$. One checks that
\[
\begin{array}{lll}
q^*_0=1,& q^*_1= (Q_{\max}\!-1)\mbox{\bf 1}_{\{V_1\ge 0\}} + (Q_{\min}\!-1) \mbox{\bf
1}_{\{V_1< 0\}}  &\mbox{ on }\; \{ \E(V_1|\F_0)\le V_0 \},\\
q^*_0= Q_{\max}\!-1, &  q^*_1 =  \mbox{\bf 1}_{\{V_1\ge 0\}} + (Q_{\min}\!-Q_{\max}+1) \mbox{\bf
1}_{\{V_1< 0\}}   &\mbox{ on }\{- \E(V^-_1|\F_0)\le V_0< \E(V_1|\F_0) \}\\
 q^*_0= Q_{\min}\!-1, &  q^*_1 =  1  &\mbox{ on }\{ V_0<- \E(V^-_1|\F_0) \}.
\end{array}
\]
Once again on the event
$$
\{ -\E(V_1^-|\F_0)  \le
V_0<\E(V_1|\F_0)\}\cap \{V_1<0\}
$$
the optimal control is not bang-bang, except if
$Q_{\max}\!\in\{1,2\}$ or $Q_{\max}=Q_{\min}$ or $Q_{\max}= Q_{\min}+1$.

\ms Finally, note that when $V_0,\, V_1\ge 0$, the events on which the optimal controls are
not bang-bang are empty. $\cqfd$

\subsection{Proof of Theorem~\ref{main}}

 $(a)$  We proceed by induction on $n$.  For $n=1$ the result is trivial
since  $T^+(1)= T^+_{00}$ and $P_0(Q)= Q_{\max}\mbox{\bf 1}_{\{V_0\ge
0\}}+ Q_{\min}\mbox{\bf 1}_{\{V_0<0\}}$. (When $n=2$ this follows from
Proposition~\ref{n=2}.)

\ss
Now, we pass from $n$ to $n+1$. Note that combining
the backward programming principle and $\mbox{\bf P1}$ yields
\begin{equation}\label{PgmDyn3}
P^{n+1}_0(Q,{\F}) = \sup\left\{xV_0+ \E\left(P^n_0(\chi^{n}(Q,x),
(\F_{1+\ell})_{0\le \ell\le n-1})\,|\,\F_0\right),\;
x\!\in I^n_Q\right\}.
\end{equation}

We inspect successively all the triangles of the tiling of $T^+(n+1)$ as follows: the upper and lower triangles
which lie strictly
inside the tiling, then the triangles which lie on the boundary of the tiling.

\bs \ni \fbox{$Q\!\in T ^+_{i\,j},\; 1\le i\le j\le n-1$:} Then,
$\chi^n(Q,x)= Q -x(1,1)$ and $I^n_Q = [0,1]$ . One checks that
$\chi^n(Q,x)\!\in T ^+_{i\,j}$ if $x\!\in[0, Q_{\min}-i]$,
$\chi^n(Q,x)\!\in  T^-_{i-1,j} $ if
$x\!\in[Q_{\min}-i,Q_{\max}\!-\!j]$  and $\chi^n(Q,x)\!\in
T^+_{i-1\,j-1}$ if $x\!\in[Q_{\max}\!-\!j,1]$ (see Figure
\ref{figure2a}). These three triangles $T ^+_{i\,j}$, $T^-_{i-1\,j}$
and $T^+_{i-1\,j-1}$ are included in $T^+(n)$. It follows from the
induction assumption that $(u,v)\mapsto P^n_0((u,v), (\F_{1+.}))$ is
$a.s.$ affine on them. Hence there exists three triplets of
$\F_1$-measurable random variables $(A^m, B^m, C^m)$, $m=1,2,3$,
such that, for every $Q\!\in T^+_{i\,j}$,
\[
P^n_0(\chi^n(Q,x),(\F_{1+.}))= \sum_{m=1}^3 \mbox{\bf
1}_{J^m_Q}(A^m(Q_{\min}-x)+B^m(Q_{\max}-x) +C^m)
\]
where $J^1_Q= [0, Q_{\min}-i]$, $J^2_Q=[Q_{\min}-i,Q_{\max}-j]$ and $J^3_Q=[Q_{\max}-j,1]$. Note these  random
coefficients satisfy
some compatibility constraints to ensure  concavity (and continuity). Consequently
\[
xV_0+ \E(P^n_0(\chi^n(Q,x),(\F_{1+.}))\,|\, \F_0) = \sum_{m=1}^3
\mbox{\bf 1}_{J^m_Q}(xV_0+ A_0^m(Q_{\min}-x)+B_0^m(Q_{\max}-x) +C_0^m)
\]
where $A^m_0= \E(A^m\,|\, \F_0)$, etc. A piecewise affine function reaches its maximum on a compact interval either
at its endpoint
or at its monotony breakpoints $x_1=0$, $x_2=Q_{\min}-i$, $x_3=Q_{\max}-j$, $x_4=1$. Hence,
\begin{eqnarray*}
\sup_{x\in I^n_Q}\left(xV_0+ \E(P^n_0(\chi^n(Q,x),(\F_{1+.}))\,|\, \F_0)\right)& =& \max \left\{x_\ell V_0+
A_0^m(Q_{\min}-x_\ell)+B_0^m(Q_{\max}-x_\ell) +C_0^m,\right.\\
&& \hskip 1 cm \left.(\ell,m)=(1,1),\,(2,1),\,(3,2),\,(4,3)\right\}.
\end{eqnarray*}

It is clear that the right hand side of the previous equality stands as the maximum of four affine functions of $Q$. One
derives that
$Q\mapsto P^{n+1}_0(Q,\F)$ is a convex function on
$T^+_{i\,j}$ as the maximum of affine functions. Hence it is affine since we know that it is also concave.

\setlength{\unitlength}{1mm}
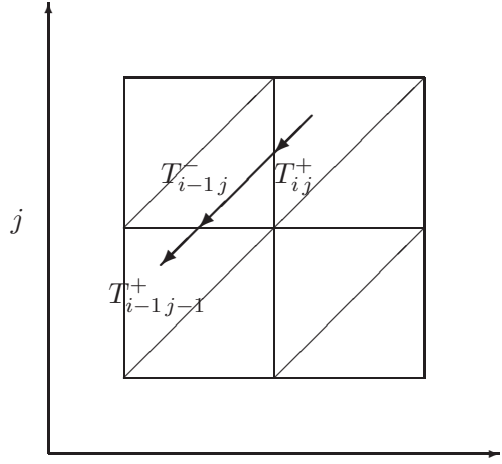
\begin{figure}
\begin{center}
\begin{picture}(60,60)
\put(0,0){\vector(0,1){60}} \put(0,0){\vector(1,0){60}}

\put(10,10){\line(0,1){40}} \put(10,10){\line(1,0){40}}
\put(10,10){\line(1,1){40}} \put(10,30){\line(1,0){40}}
\put(30,10){\line(0,1){40}} \put(10,50){\line(1,0){40}}
\put(50,10){\line(0,1){40}} \put(30,10){\line(1,1){20}}
\put(10,30){\line(1,1){20}} \linethickness{2.5mm} \thicklines
\put(20,30){\vector(-1,-1){5}} \put(30,40){\vector(-1,-1){10}}
\put(35,45){\vector(-1,-1){5}}

\put(0,-5){$0$} \put(-5,30){$j$} \put(30,-5){$i$}
\put(30,36){$T^{+}_{i\,j}$} \put(15,36){$T ^{-}_{i-1\,j}$}
\put(8,20){$T ^{+}_{i-1\,j-1}$}

\end{picture}
\caption{$x\mapsto\chi^n(x)\text{ for }Q\!\in T ^+_{i\,j},\; 1\le
i\le j\le n-1$} \label{figure2a}
\end{center}
\end{figure}

\bs
\ni \fbox{$Q\!\in T ^-_{i\,j},\; 1\le i<j\le n-1$:} This case can be treated
likewise.

\bs
\ni \fbox{$Q\!\in T ^{\pm}_{0\,j},\,1\le j\le n-1$:} In that case $I^n_Q=
[0,1]$, $\chi^n(Q,x)= ((Q_{\min}-x)^+,
Q_{\max}-x),
\, x\!\in I^n_Q$.

\ss -- If $Q\!\in T^+_{0\,j}$, $\chi^n(Q,x)=Q-x(1,1)\!\in T^+_{0j}$,
$x\!\in[0,Q_{\min}]$, $\chi^n(Q,x)=(0,Q_{\max}-x)\!\in T^+_{0\,j}$,
$x\!\in [Q_{\min},Q_{\max}-j]$, $\chi^n(Q,x)=(0,Q_{\max}-x)\!\in
T^+_{0\,j-1}$, $x\!\in [Q_{\max}-j,1]$ (see Figure \ref{figure2b}).
The induction assumption implies that $x\mapsto
P^n_0(\chi^n(Q,x),(\F_{1+.}))$ is piecewise affine with monotony
breaks at $Q_{\min}$ and $Q_{\max}-j$.

\ss --  If $Q\!\in T^-_{0\,j}$, $\chi^n(Q,x)$ crosses the upper
(horizontal) edge of $T^+_{0\,j-1}$ at $x=Q_{\max}-j$ and the left
(vertical) edge of $T^+(n)$ at $x=Q_{\min}$. Hence $x\mapsto
P^n_0(\chi^n(Q,x),(\F_{1+.}))$ is again piecewise affine with
monotony breaks at $Q_{\min}$ and  $Q_{\max}-j$.

\ss
In both cases one concludes as above.

\setlength{\unitlength}{1mm}
\begin{figure}
\begin{center}
\begin{picture}(60,60)
\put(0,0){\vector(0,1){60}} \put(0,0){\vector(1,0){50}}

\put(0,30){\line(1,0){20}} \put(0,50){\line(1,0){40}}
\put(0,10){\line(1,1){40}} \put(0,30){\line(1,1){20}}
\put(20,30){\line(0,1){20}} \linethickness{2.5mm} \thicklines
\put(0,30){\vector(0,-1){10}} \put(0,40){\vector(0,-1){10}}
\put(5,45){\vector(-1,-1){5}}

\put(0,-5){$0$} \put(-5,30){$j$} \put(8,26){$T ^{+}_{0\,j-1}$}
\put(8,46){$T ^{+}_{0\,j}$}

\end{picture}
\caption{$x\mapsto\chi^n(x)\text{ for }Q\!\in T ^{+}_{0\,j},\,1\le
j\le n-1$} \label{figure2b}
\end{center}
\end{figure}

\bs
\ni \fbox{$Q\!\in T ^{\pm}_{0\,0}$:} On proceeds like with $T^+_{0j}$
except that $I^n_Q=
[0,Q_{\max}]$ which yields only one monotony break at $Q_{\min}$.

\bs \ni \fbox{$Q\!\in T ^{\pm}_{i\,n},\,1\le i\le n-1$:} Assume
first $Q\!\in T ^{+}_{i\,n}$. $I^n_Q=[0,1]$ and $\chi^n(Q,x)=
(Q_{\min}-x,n)$ if $x\!\in [0, Q_{\max}-n]$.  Otherwise
$\chi^n(Q,x)= (Q_{\min}-x,Q_{\max}-x)$. It follows (see Figure
\ref{figure3}) that $\chi^n(Q,x)\!\in T^+_{i,n-1}$ if $x\!\in
[0,Q_{\min}-i]$ and $\chi^n(Q,x)\!\in T^+_{i-1,n-1}$ if $x\!\in
[Q_{\min}-i,1]$. Both $T^+_{i\,n-1}$ and $T^+_{i-1\,n-1}$ are
included in $T^+(n)$. Hence $(u,v)\mapsto P^n_0((u,v),(\F_{1+.}))$
is   affine on  both triangles, one derives that
\begin{eqnarray*}
 P^n_0\!\!\!\!\!&\!\!\!(\!\!\!&\!\!\!\!\!\!\chi^n(Q,x),(\F_{1+.}))\\
&=& \mbox{\bf 1}_{_{x\in[0,
Q_{\max}-n]}} (A^1(Q_{\min}-x)+B^1n+C^1)+\mbox{\bf 1}_{_{x\in[Q_{\max}-n,1]}}(A^2(Q_{\min}-x)+B^2(Q_{\max}-x)+C^2).
\end{eqnarray*}
where $A^m,B^m,C^m,\, m=1,2$ are $\F_1$-measurable r.v.. Then, one concludes like in the first case.

If $Q\!\in T ^{-}_{i\,n}$, one proceeds likewise except that the two ``visited" triangles of $T^+(n)$ are
$T^{\pm}_{i-1,n-1}$.

\setlength{\unitlength}{1mm}
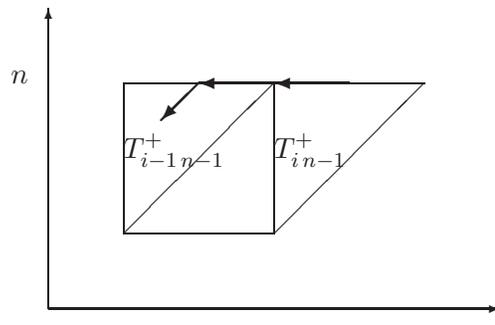
\begin{figure}
\begin{center}
\begin{picture}(60,60)
\put(0,0){\vector(0,1){40}} \put(0,0){\vector(1,0){60}}

\put(10,10){\line(0,1){20}} \put(10,10){\line(1,0){20}}
\put(10,10){\line(1,1){20}} \put(10,30){\line(1,0){40}}
\put(30,10){\line(1,1){20}} \put(30,10){\line(0,1){20}}
\linethickness{2.5mm} \thicklines \put(20,30){\vector(-1,-1){5}}
\put(30,30){\vector(-1,0){10}} \put(40,30){\vector(-1,0){10}}

\put(0,-5){$0$} \put(-5,30){$n$} \put(30,-5){$i$} \put(30,20){$T
^{+}_{i\,n-1}$} \put(10,20){$T ^{+}_{i-1\,n-1}$}
\end{picture}
\caption{$x\mapsto\chi^n(x)\text{ for }Q\!\in T ^{+}_{i\,n},\,1\le
i\le n-1$} \label{figure3}
\end{center}
\end{figure}

\bs
\ni \fbox{$Q\!\in T ^{\pm}_{0\,n}$:} $I^n_Q=[0,1]$ and

\ss
-- $\chi^n(Q,x)= ((Q_{\min}-x)^+,n), \, x\!\in I^n_Q$ if $Q\!\in T
^{+}_{0\,n}$,

\ss -- $\chi^n(Q,x)= (Q_{\min}-x,n), \, x\!\in [0,Q_{\max}-n]$,
$\chi^n(Q,x)= (Q_{\min}-x,Q_{\max}-x), \, x\!\in
[Q_{\max}-n,Q_{\min}]$ if $Q\!\in T^{-}_{0\,n}$.

\ss
In both cases the only ``visited" triangle is $T^+_{0\,n-1}\subset
T^+(n)$ and one concludes as usual.

\bs
\ni \fbox{$Q\!\in T ^{+}_{nn}$:} $I^n_Q=[Q_{\min}-n,1]$ and $\chi^n(Q,x)=
(Q_{\min}-x,n)$ if $x\!\in [Q_{\min}-n,Q_{\max}-n]$, $\chi^n(Q,x)=
(Q_{\min}-x,Q_{\max}-x)$ otherwise. Hence $\chi^n(Q,x)$ takes its values
in $T^+_{n-1\,n-1}$ on which $(u,v)\mapsto P^n_0((u,v),(\F_{1+.}))$ is
affine. The conclusion follows.

The inspection of all these cases completes the proof of the induction.

\ms The values of $P^n_0(Q)$ when $Q\!\in\{(0,0),\, (0,n),\, (n,n)\}$
are obvious consequences of the degeneracy of the global constraints.

\bs
\ni $(b)$  We deal successively
with the two announced settings.

\ms
\ni -- {\sc global constraints in $\N^2$:} Let $n\ge1$. We rely on the
characterization~(\ref{Optiqk}) of $q^*_0$
\[
q^*_0= {\rm argmax}_{x\in I^{n-1}_Q}\left(xV_0+ \E\left(P^{n-1}_0(\chi^{n-1}(Q,x), (\F_{1+\ell})_{0\le \ell\le
n-2})\,|\,\F_0\right)\right).
\]
We know from item~$(a)$ that $(u,v)\mapsto P^{n-1}((u,v),
\F_{1+.})$ is    affine on every tile $T^\pm_{ij}$ of  $T^+(n)$.

If $Q\!\in T^+(n)\cap \N^2$ then one checks that
$x\mapsto \chi^{n-1}(Q,x)$, $x\!\in I^{n-1}_Q$,  is always affine with $I^{n-1}_Q$ having  $0$
and/or $1$ as endpoints. To be precise

\ms
-- if $Q=(i,j)$, $1\le i\le j\le n-1$, $\chi^{n-1}(Q,x)= (i-x,j-x)\!\in\partial T^+_{i-1\,j-1} \cap \partial T^-_{i-1\,j-1}$,
$I^{n-1}_Q=[0,1]$,

\ms
-- if $Q=(0,j)$, $1\le j\le n-1$, $\chi^{n-1}(Q,x)= (0,j-x)\!\in \partial T^+_{0\,j-1}$, $I^{n-1}_Q=[0,1]$,

\ms
-- if $Q=(0,0)$,   $\chi^{n-1}(Q,x)= (0,0)$, $I^{n-1}_Q=\{0\}$,

\ms
--  if $Q=(i,n)$, $1\le i\le   n-1$,  $\chi^{n-1}(Q,x)= (i-x,n-1)\!\in \partial T^-_{i-1\, n-1} $, $I^{n-1}_Q=[0,1]$,

\ms
--  if $Q=(0,n)$,    $\chi^{n-1}(Q,x)= (0,n-1)$, $I^{n-1}_Q=[0,1]$,

\ms
--  if $Q=(n,n)$,    $\chi^{n-1}(Q,x)= (n-1,n-1)$, $I^{n-1}_Q=\{1\}$.

\ms As a consequence,  affinity being stable by composition,
$x\mapsto P_0^{n-1}(\chi^{n-1}(Q, x),\F_{1+.})$ is affine on
$I^{n-1}_Q\in \{[0,1], \{0\}, \{1\}\}$. In turn,

$x\mapsto xV_0+\E(P_0^{n-1}(\chi^{n-1}(Q, x),\F_{1+.})\,|\,\F_0)$
is affine and   reaches its maximum at some endpoint of
$I^{n-1}_Q$ $i.e.$   $q^*_0=0$ or at $q^*_0=1$. Then, inspecting
the above cases shows that $Q^{1,*}=Q-q^*_0(1,1)\!\in T^+(n-1)\cap
\N^2$. Using~(\ref{PgmDynGlob}) and~(\ref{Optiqk}), one shows  by
induction on $k$ that $q^*_k$ is always $\{0,1\}$-valued.

 \ms \ni
-- {\sc non negative $V_i$:}

{\em Step~1: Global constraint saturated.} Let $n\!\ge 1$. Let
$q^*=(q^*_k)_{0\le k\le n-1}$ be an optimal $Q$-admissible
control. We introduce the $\F$-stopping time
\[
\tau(q^*) := \min\left\{k\,|\, q^*_0+\cdots+q^*_k <Q_{\max}-(n-1)+k\right\}
\]
with the convention $\min\emptyset =+\infty$.

On  $\tau(q^*)=+\infty$, $q^*_0+\cdots+q^*_k\ge Q_{\max}-(n-1)+k$ for every $k=0,\ldots,n-1$.
In particular the global constraint is saturated at time $n-1$, $i.e.$
\[
q^*_0+\cdots+q^*_{n-1}=Q_{\max}.
\]
Set
\[
\tilde q_k=q^*_k\mbox{\bf 1}_{\{k\leq \tau(q^*)-1\}}+(
Q_{\max}-(n-1)+\tau(q^*)-(q^*_0+\cdots+q^*_{\tau(q^*)-1}))\mbox{\bf 1}_{\{k=\tau(q^*)\}} .
\]

One check that   $\tilde q$ is a $Q$-admissible control: this
follows from the fact that $\tau(q^*)$ is an $\F$-stopping time
(note that $q^*_\tau$ is $[0,1]$-valued on $\{\tau(q^*)
<+\infty\}$ since  then $q^*_0+\cdots+q^*_{\tau-1} \ge
Q_{\max}-n+\tau$ and $q^*_0+\cdots+q^*_\tau<Q_{\max}-(n-1)+\tau$).
Likewise one shows that $\tilde q_k\ge q^*_k$ for every
$k=0,\ldots,n-1$.

Furthermore, note that if $Q_{\max}$ is an integer and $q^*$ is $\{0,1\}$-valued then $\tilde q$ is still
$\{0,1\}$-valued.

 The
$V_k$ being non negative
\[
\sum_{k=0}^{n-1} \tilde q_kV_k\ge \sum_{k=0}^{n-1}  q^*_kV_k
\]
hence $\tilde q$ is still an optimal control. Furthermore $\tilde q_0+\cdots+\tilde q_k \ge Q_{\max}-(n-1)+k$ on
$\{k\ge \tau(q^*)\}$  so that  the stopping time $\tau(\tilde q)$ satisfies by construction
\[
\tau(\tilde q) \le \tau(q^*)-1 \;\mbox{ on }\; \{1\le
\tau(q^*)<+\infty\}\quad\mbox{ and }\quad \tau(\tilde q) =+\infty
\; \mbox{ on }\;\{\tau(q^*)=0\}\cup\{\tau(q^*)=+\infty\}.
\]
Note that if the control $q^*$ is bang-bang, iterating the above construction  at most $n$ times yields
an optimal control $q^{opt}$ for
which
$\tau(q^{opt})=+\infty$
$a.s.$. Such a control $q^{opt}$ saturates the global constraint.

As a consequence, this shows that $P_0(Q,\F)=P_0((0,
Q_{\max}),\F)= P_0((Q_{\max}, Q_{\max}),\F)$ $a.s.$ so that
$Q_{\max}\mapsto P_0((Q_{\max}, Q_{\max}),\F)$ is $a.s.$
non-decreasing and concave.

\ms
{\em Step~2: Local constraints.} Since there is an optimal control $q^*$
which saturates the global constraint, one may assume without loss of
generality that $Q_{\min}=Q_{\max}$. We proceed again by induction on $n$ based
on the dynamic programming formula~(\ref{PgmDyn3}). When
$n=1$ the result is obvious (and true when $n=2$ as well).

Assume now  the
announced result is true for $n\ge1$.

Let $j\!\in\{0,\ldots, n-1\}$ and $Q_{\max}\!\in [j,j+1]$. Then,
$I^n_Q=[0,Q_{\max}\wedge\!1]$ and $\chi^n((Q_{\max},Q_{\max}),x)
=(Q_{\max}-x,Q_{\max}-x)$, $x\!\in I^n_Q$. Hence
$\chi^n((Q_{\max},Q_{\max}),x)\!\in T^+_{jj}$, $x\!\in
[0,Q_{\max}-j]$ and $\chi^n((Q_{\max},Q_{\max}),x)\!\in
T^+_{j-1,j-1}$, $x\!\in [Q_{\max}-j,Q_{\max}\!\wedge\!1]$. Now
$v\mapsto P^n_0((v,v),\F_{1+.})$ is $a.s.$ concave,
non-decreasing, affine on $[j-1,j]$ and on $[j,j+1]$ and
non-decreasing. Consequently, there exists $B^m,\, C^m$, $m=1,2$,
$\F_1$-measurable random variables satisfying
\begin{eqnarray*}
P_0(\chi^n((Q_{\max},Q_{\max}),x),\F_{1+.}) &=& B^1(Q_{\max}-x)+C^1,\;   x\!\in
[0,Q_{\max}-j],\\ P_0(\chi^n((Q_{\max},Q_{\max}),x),\F_{1+.})& =& B^2(Q_{\max}-x)+C^2,\;
  x\!\in [Q_{\max}-j,Q_{\max}\!\wedge\!1].
\end{eqnarray*}
with $0\le B^1\le B^2$ and $ B^2j+C^2= B^1j+C^1$ $a.s.$. Set temporarily
$$\Psi(x) :=
xV_0+\E\left(P_0(\chi^n((Q_{\max},Q_{\max}),x),\F_{1+.})\,|\,
\F_0\right).
$$
 Hence,
\[
\sup_{x\in I^n_Q}\Psi(x)=\max(\Psi(0), \Psi(Q_{\max}-j),\Psi(Q_{\max}\wedge 1)).
\]

Set $B^m_0:= \E(B^m\,|\,\F_0)$ and $C^m_0:=\E(C^m\,|\,\F_0) $ and note that $B^1_0\le
B^2_0$ and  $ B_0^2\,j+C_0^2= B_0^1\,j+C_0^1$ $a.s.$. Elementary computations show that:

\ms
--  $\Psi(0)\le \Psi(Q_{\max}-j)$ on $\{V_0\ge B^1_0 \}$  and $\Psi(0)\ge
\Psi(Q_{\max}-j)$ on $\{V_0\le B^1_0 \}$,

\ms
--  $\Psi(Q_{\max}-j)\le \Psi(Q_{\max}\wedge 1)$ on $\{V_0\ge B^2_0 \}$  and $\Psi(Q_{\max}-j)\ge
\Psi(Q_{\max}\wedge 1)$ on $\{V_0\le B^2_0 \}$.

\ms
Consequently $q^*_0$ can be chosen   $\{0,Q_{\max}\wedge 1\}$-valued on
$E_0:=\{B^2_0<V_0\le  B^1_0\}\!\in
\F_0$ and equal to $Q_{\max}-j$ on $^cE_0:=\{ B^1_0 <V_0\le B^2_0\}\!\in \F_0$.

\ms On $E^1_0=E_0\cap\{q^*_0=Q_{\max}\wedge 1\}\!\in \F_0$, one
has $P^{n+1}_0((Q_{\max}, Q_{\max}),\F)=P^{n+1}_0((Q_{\max},
Q_{\max}),\F\cap E^1_0)$.

Then, the  dynamic programming formula shows that the other components $(q^*_k)_{1\le
k\le n}$ of the optimal control on $E_0$  can be obtained as the optimal control of
the pricing problem $P_0^{n}( ((Q_{\max}\!-\!1)^+,(Q_{\max}\!-\!1)^+),(\F_{1+k}\cap E^1_0)_{0\le k\le
n-1})$. One derives from the induction assumption at time $n$ that $(q^*_k)_{1\le
k\le n}$ can be chosen  bang-bang and $((Q_{\max}\!-\!1)^+,(Q_{\max}\!-\!1)^+)$-admissible which implies that $q$ is
$Q$-admissible and bang-bang since $q^*_0=Q_{\max}\wedge 1$ (on $E^1_0$). A similar proof holds
on $E^0_0=E_0\cap\{q^*_0=0\}$.

\ms
On $^c E_0$, one has $P^{n+1}_0((Q_{\max}, Q_{\max}),\F)=P^{n+1}_0((Q_{\max},
Q_{\max}),\F\cap ^c E_0)$. Then, the  dynamic programming formula shows that the other components
$(q^*_{1+k})_{0\le k\le n-1}$ of the optimal control on $^cE_0$  can be obtained as the optimal control of
the pricing problem $P_0^{n}( (j,j),(\F_{1+k}\cap ^cE_0)_{0\le k\le
n-1})$. As $(j,j)\!\in \N^2$ there exists a $(j,j)$-admissible bang-bang optimal control
$(q^*_k)_{0\le k\le n-1}$ (with respect to $(\F_{1+k}\cap ^cE_0)_{0\le k\le
n-1})$  on $^c E_0$. Then $q^*_{1+k}$ is $\{0,1\}$-valued for every $k=0,\ldots, n-1$ (in fact
identically $0$ if $j=0$). At this
stage one can recursively modify $(q_{1+k}^*)_{0\le k\le n-1}$  using the procedure described in Step~1 to
saturate  the upper
global constraint. Finally one may assume that $\sum_{0\le k\le n-1}q_k^*=j$ which in turn implies that
$(q^*_k)_{0\le k\le n}$ is a bang-bang
$(Q_{\max}, Q_{\max})$-admissible optimal control.$\quad\cqfd$

\bs \ni {\sc Application.} When a global constraint $Q$ belongs to
the interior of a triangle $T^\pm_{i\,j}$,  {\em one only needs to
compute the value of $P_0(.,\F)$ at the vertices of this triangle to
derive the value of the premium at every  $Q\!\in T^\pm_{i\,j}$}.
When $Q$ is itself an integral valued couple, at most six further
points allow to compute the premium in a neighborhood of $Q$. We
will use this result extensively when designing our quantization
based numerical procedure in Section~\ref{SwingQuant}.

\bs \noindent{\sc An additional result.}  Proposition~\ref{n=2}
shows that it is hopeless to produce in full generality some
bang-bang optimal control when $Q_{\max}-Q_{\min}\!\in\N$. This
comes from the fact that at integral valued global constraints the
bang-bang  optimal control may saturate none of the global
constraints (indeed, so is the case at $(2,2)$ when $n=2$).
However, using the same approach as that developed in that in the
proof of Theorem~\ref{main}, one can show the following
result, whose details of proof are left to the reader.

\begin{Cor} Assume the assumptions of Theorem~\ref{main} hold. If a
couple of admissible constraint $(Q_{\min},Q_{\max})\!\in T^+(n)$
satisfies
$$
Q_{\max}-Q_{\min}\!\in \{0,\ldots,n\}
$$
then there exist a {\em quasi-bang-bang} control in the following
sense: $\P$-$a.s.$, $q^*_k$ is $\{0,1\}$-valued except for at most
one local constraint $q^*_{k_0}$.
\end{Cor}


\subsection{The Markov setting}

\setcounter{equation}{0}
\setcounter{Assumption}{0}
\setcounter{Theorem}{0}
\setcounter{Proposition}{0}
\setcounter{Corollary}{0}
\setcounter{Lemma}{0}
\setcounter{Definition}{0}
\setcounter{Remark}{0}
By Markov setting we simply mean that the payoffs $V_k$ are
function of an $\R^d$-valued underlying $\F$-Markov structure
process $(Y_k)_{0\le   k\le n-1}$ $i.e.$
\[
V_k= v_k(Y_k),\; k=0,\ldots,n-1.
\]
The Markovian dynamics of $Y$ reads on Borel functions $g:\R^\to\R$
\[
\E(g(Y_{k+1})\,|\, \F_k)=\E(g(Y_{k+1})\,|\, Y_k)= \Theta_k(g)(Y_k)
\]
where $(\Theta_k)_{0\le k\le n-1}$ is  sequence of Borel probability transitions on $(\R^d, \B(\R^d))$.

 Then the backward dynamic programming principle~(\ref{PgmDyn}) can be rewritten  as follows
\[
P^n_k(Q)= p^n_k(Q,Y_k),\; k=0,\ldots,n
\]
with $p^n_n(.,.) \equiv 0$  and for every $k=0,\ldots,n-1$ and
every  $Q\!\in T^+(n\!-\!k)$,
\begin{eqnarray}\label{petitPGD}
p^n_k(Q,y)&= & \sup \left\{x\,v_{k}(y) +
\Theta_{k}(p^n_{k+1}(\chi^{n\!-\!k-1}(Q,x),.))(y),\, x\!\in  I^{n\!-\!k-1}_Q\right\},
\end{eqnarray}
where $\chi^{M}(Q,x)= ((Q_{\min}-x)^+,
(Q_{\max}-x)\!\wedge\!M)$ and $I^M_Q:=[(Q_{\min}-M)^+\!\wedge\! 1,Q_{\max}\!\wedge\! 1]$.

\ms
\noindent {\bf Pointwise estimation of $P^n_0(Q^0)$.} As established in Theorem~\ref{main}, one only needs to
compute the value function
$P^n_0(Q)$ at global constraints $Q=(Q_{\min},Q_{\max})\in T^+(n)\cap \N^2 $ $i.e.$ with integral components.
Moreover, for these
constraints, the local optimal control $q^*_k$ is always bang-bang $i.e.$ $q^*_k\!\in\{0,1\}$.

Let $Q^0=(Q^0_{\min}, Q^0_{\max})\!\in T^+(n)\cap \N^2$. For every
$k=0,\ldots,n-1$, one defines the set of attainable residual global
constraints at time $k$, namely
\begin{equation}\label{reachable}
{\cal Q}^n_k(Q^0):=\left\{((Q^0_{\min}-\ell)^+, (Q^0_{\max}-\ell)^+\!\wedge\!(n\!-\!k)),\; \ell=0,\ldots,k\right\}.
\end{equation}
(thus ${\cal Q}^n_0(Q^0)=\{Q^0\}$). Note that the running
parameter $\ell$ represents the possible values of the cumulative
purchase process $\bar q^*_k$.

One checks that for every $Q\!\in {\cal Q}^n_k(Q^0)$,
$\chi^{n\!-\!k-1}(Q,1)\chi^{n\!-\!k-1}(Q,0)\!\in {\cal
Q}^n_{k+1}(Q^0)$ since
$$
\chi^{n\!-\!k-1}(Q,1)=((Q_{\min}-1)^+, (Q_{\max}-1)^+)\;\mbox{ and
}\;\chi^{n\!-\!k-1}(Q,0)=(Q_{\min},Q_{\max}\wedge(n\!-\!k-1)).
$$
Consequently the backward
dynamic programming formula having $p^n_0(Q^0,y)$ as a result reads:
\[
 p^n_k(Q,y)=\max \left\{x\,v_{k}(y) +
\Theta_{k}(p^n_{k+1}(\chi^{n\!-\!k-1}(Q,x),.))(y),\,x\!\in\{0,1\}\cap I_Q^{n-1-k}\!\right\},Q\!\in{\cal Q}^n_k(Q^0),\, k=0,\ldots,n-1.
\]

At this stage no numerical computation is possible yet  since   no space discretization has been achieved. This is
the aim of
the  next section where we will approximate the above dynamic programming principle by (optimal) quantization of the
state process
$Y$.

\section{Computing swing contracts by (optimal) quantization}\label{SwingQuant}

\subsection{The abstract quantization tree approach}
\paragraph{Abstract quantization} In this section, we propose a quantization based approach to compute the premium
of the swing contracts
with firm constraints. Quantization has been originally introduced and developed in the early 1950'  in Signal
processing (see~\cite{IEEE}).
 The starting idea is simply to replace every random vector $Y:(\Omega,{\cal A})\to \R^d$ by a random vector $
 \widehat Y= g(Y)$ taking finitely many values in a {\em grid} (or {\em codebook}) $\Gamma:=\{y_1,\ldots, y_{N}\}$
 (of size $N$). The  grid $\Gamma$ is also called an $N$-quantizer of $Y$. When the Borel  function $g$ satisfies
\begin{equation}
|Y-\widehat Y| =d(Y, \Gamma)=\min_{1\le i\le N}|Y-y^i|,
\end{equation}
$\widehat Y$ is called a {\em Voronoi quantization} of $Y$  (and $g$
as well). One easily checks that $g$ is necessarily a nearest
neighbor projection on $\Gamma$ $i.e.$ satisfies
\[
\forall\, i\!\in\{1,\ldots, N\}, \qquad \{g=y^i\}\subset \{u\!\in
\R^d\,:\, |u-y^i|=\min_{1\le j\le N}|u-y^j|\}.
\]
 The so-called {\em Voronoi cells} $\{g=y^i\}),\,1\le i\le N$, make up a {\em Voronoi tessellation} or partition of
 $\R^d$ (induced by $\Gamma$. Note that when the distribution $\P_{_Y}$ of $Y$ weights no
hyperplanes the boundary of the Voronoi tessellation of $\Gamma$
are $\P_{_Y}$-negligible so that the $\P_{_Y}$-weights of the
Voronoi cells entirely characterize the distribution of $\widehat
Y^\Gamma$.

 When $p\!\in [1,\infty)$, the $L^p$-mean error induced by replacing $Y$ by $\widehat Y$, namely
$$
\|Y-\widehat Y\|_{_p}= \left(\E(\min_{1\le i\le
N}|Y-y^i|^p\right)^{\frac 1p}
$$
is called the {\em $L^p$-mean quantization error} induced  by $\Gamma_k$ and its $p^{th}$ power is known as the
$L^p$-distortion.  We will see in the next section that the codebook $\Gamma$ can be optimized so as to minimize the
$L^p$-quantization error with respect to the (distribution of) $Y$.

\medskip
Our aim in this section is to design an algorithm based on the
quantization of the Markov chain $(Y_k)$ at every time $k$ to
approximate the premium of the swing contract with firm
constraints and to provide some  {\em a priori} error estimates in
terms of quantization errors.

\paragraph{Quantized tree for pricing swing options.} As a first step we consider at every time $k$ a grid
$\Gamma_k:=\{y_k^1,\ldots,y_k^{N_k}\}$ (of size $N_k$). Then, we
design the quantized tree algorithm to price swing contracts by
simply mimicking the original dynamic programming
formula~(\ref{PgmDyn}). This means in particular that we force in
some way the Markov property on $(\widehat Y_k)_{0\le k\le n-1}$
by considering the quantized transition operator
\[
\widehat \Theta_k(g)(y_k^i)=
\sum_{j=1}^{N_{k+1}}g(y_{k+1}^j)p^{ij}_{k},\qquad p^{ij}_{k}
:=\P(\widehat Y_{k+1}=y_{k+1}^j\,|\, \widehat Y_k=y_k^i)
\]
so that
\[
\widehat \Theta_k(g)(\widehat Y_k) = \E\left(g(\widehat Y_{k+1})\,|\,\widehat Y_{k}\right), k=0,\ldots,n-1.
\]
$\rhd$ Let $Q^0\!\in T^+(n)\cap \N^2$ be a couple of (deterministic)
global constraints (at time $0$).  The quantized dynamic
programming principle  is  defined by
\begin{eqnarray}\label{PgmDynQuantGlob0}
\nonumber \widehat P^n_n(Q)\!\!&\!\!:=\!\!&\!\!0,\quad Q\!\in T^+(n)\cap \N^2,\\
 \hskip -0.35 cm \widehat P^n_k(Q) \!\!&\!\!:=\!\!&\!\!\max\!\left(\!x\,v_k(\widehat Y_k) \!+\!
\E\!\left(\!\widehat P^n_{k+1}\!(\chi^{n\!-\!k-1}\!(Q,x))\,|\,\widehat Y_k\!\right)\!,
x\!\in\partial I^{n-1-k}_Q\!\right)\!,Q\!\in{\cal Q}^n_k(Q^0),  k=0,\ldots,n-1.\qquad
\end{eqnarray}
One easily shows by induction that, for every $ Q\!\in{\cal Q}^n_k(Q^0)$ (residual global constraint at time $k$),
\[
\widehat P^n_k(Q)= \widehat p^n_k(Q,\widehat Y_k)
\]
\begin{eqnarray}\label{PgmDynQuantGlob}
\nonumber \mbox{where }\hskip 1 cm \widehat p^n_{n}(Q,y)&=&0,\quad Q\!\in{\cal Q}^n_n(Q^0),\; y\!\in\R^d, \\
\hskip -0.35 cm \widehat p^n_k(Q,y_k^i) &= &\sup \left\{x\,v_k(y^i_k) +
\widehat \Theta_k(\widehat p^n_{k+1}(\chi^{n\!-\!k-1}(Q,x),.))(y_k^i),x\!\in\partial I^{n-1-k}_Q\right\}\hskip 1,75 cm\\
\nonumber && i=1,\ldots, N_k,\; \,Q\!\in{\cal Q}^n_k(Q^0),\;k=0,\ldots,n-1.
\end{eqnarray}

\noindent $\rhd$ When $Q^0\!\in T^+(n)\setminus \N^2$, one defines $\widehat P^n_0(Q)$ (and $\widehat p^n_0(Q,.))$ by affinity on each elementary triangle $T^\pm_{ij}$ that tiles $T^+(n)$.

\paragraph{Complexity} Let us briefly discuss the complexity of this quantized backward dynamic procedure.  Let $
k\!\in \{0,\ldots,n-1\}$. At every ``nod" $y_k^i$ the computation
of $\widehat \Theta_k(\widehat
p^n_{k+1}(\chi^{n\!-\!k-1}(Q,x),.))(y_k^i)$ requires $N_{k+1}$
products (up to a constant), so that for a given residual global
constraint the complexity at time $k$ in the dynamic programming
is proportional to $N_kN_{k+1}$.  On the other hand, one checks
that
$$
{\rm card}({\cal Q}_k^n(Q^0))=  (Q^0_{\max}\!\wedge \!k)\!+\! 1-(Q^0_{\max}-Q^0_{\min}-(n-k)-1)^+.
$$
Consequently, the  complexity of the computation of $\widehat p^n_0(Q^0,\widehat Y_0)$ is proportional to
\[
\sum_{k=0}^{n-1} {\rm card}({\cal Q}_k^n(Q^0))N_kN_{k+1}.
\]
A simple upper-bound is provided by
\[
\sum_{k=0}^{n-1}((Q^0_{\max}\!\wedge \!k)\!+\! 1)N_k\,N_{k+1}
\]
and a uniform one by
\[
\sum_{k=0}^{n-1}(k+1)N_k\,N_{k+1}.
\]
Note that this last upper bound corresponds to the complexity of the quantized version of the algorithm based on some {\em penalized global volume constraints} (see the companion paper~\cite{BABOPA1}).
\paragraph{A priori error bounds for the quantized procedure}

\begin{Thm}\label{Quantrate} Assume that the Markov process $(Y_k)_{0\le k\le n-1}$ is Lipschitz Feller
in the following sense: for every
bounded Lipschitz continuous function $g:\R^d\to \R$ and every $k\!\in\{0,\ldots,n-1\}$, $\Theta_k(g)$
is a Lipschitz function satisfying
$[\Theta_k(g)]_{\rm Lip}\le[\Theta_k ]_{\rm Lip} [g]_{\rm Lip}$. Assume that every function $v_k:\R^d\to \R$
is Lipschitz continuous
with Lipschitz coefficient $[v_k]_{\rm Lip}$. Let  $p\!\in [1,\infty)$ such that $\max_{0\le k\le n-1}|Y_k|\!\in L^p(\P)$.
Then, there exists a real constant $C_p>0$ such that
\begin{equation}\label{ratequantswing}
\|\sup_{Q\in  T_{\N}^+(n)}|\widehat P^n_0(Q) -  P^n_0(Q)|\|_{_p}\le C_{p}\sum_{k=0}^{n-1} \|Y_k-\widehat Y_k\|_{_p}
\end{equation}
\end{Thm}

\noindent {\bf Remark.} In most situations ${\cal F}_0=\{\emptyset, \Omega\}$ so that the error term
$|\widehat P^n_0(Q) -  P^n_0(Q)|$ is
deterministic.  When ${\cal F}_0$ is not trivial, it is straightforward from~(\ref{ratequantswing}) (with $p=1$) that
\[
\sup_{Q\in T^+(n)}|\E(\widehat P^n_0(Q) )- \E( P^n_0(Q))|\le C_{1}\sum_{k=0}^{n-1} \|Y_k-\widehat Y_k\|_{_1}.
\]

\bigskip We first need a    lemma about the Lipschitz regularity of the $p^n_k$ functions.

\begin{Lemma} For every $k\!\in \{0,\ldots,n-1\}$, the function  $y\mapsto p^n_k(Q,y)$ is Lipschitz on $\R^d$,
uniformly with respect to
$Q\!\in T^+(n-k)$ and its Lipschitz coefficient $[p^n_k]_{{\rm Lip},y}:=\sup_{Q\in T^+(n-k)}[p^n_k(Q,.)]_{\rm Lip}$ satisfies for every $k\!\in \{0,\ldots, n-1\}$,
\[
[p^n_{n-1}]_{{\rm Lip},y}\le [v_{n-1}]_{\rm Lip},\qquad
[p^n_k]_{{\rm Lip},y}\le [v_k]_{\rm Lip}+ [\Theta_k]_{\rm
Lip}[p^n_{k+1}]_{{\rm Lip},y}.
\]
\end{Lemma}

\ni{\bf Proof.} This follows easily by a backward induction on $k$, based on the dynamic programming
formula~(\ref{petitPGD}) and the elementary inequality
$|\sup_{i\in I}a_i-\sup_{i\in I}b_i|\le \sup_{i\in I} |a_i-b_i|$. $\cqfd$

\bs
\noindent {\bf Proof of Theorem~\ref{Quantrate}.}  First note that, by piecewise affinity of  $\widehat P^n_0$  and $ P^n_0$, one has
\[
\sup_{Q\in  T_{\N}^+(n)}|\widehat P^n_0(Q) -  P^n_0(Q)|= \sup_{Q\in  T_{\N}^+(n)\cap \N^2}|\widehat P^n_0(Q) -  P^n_0(Q)|.
\]

Temporarily set $T^+_\N(n):=T^+(n)\cap \N$. Let $k\!\in\{0,\ldots,n-1\}$. Now
\begin{eqnarray}
\sup_{Q\in T_\N^+(n-k)}|p^n_k(Q,Y_k)\!\!&\!\!-\!\!&\!\!\widehat p^n_k(Q,\widehat Y_k)|\le |v_k(Y_k)-v_k(\widehat Y_k)| \\
\nonumber &+& \hskip -1 cm\sup_{Q\in T_\N^+(n-k), \,x \in \partial I^{n-1-k}_Q}\left|\E(p^n_{k+1}(\chi^{n-1-k}(Q,x),Y_{k+1})\,|\F_k)-\E(\widehat  p^n_{k+1}(\chi^{n-1-k}(Q,x),\widehat Y_{k+1})\,|\, \widehat Y_k) \right|.
\end{eqnarray}

Now, using that $\Theta_k$ is a Markov transition and that $\sigma(\widehat Y_k)\subset \sigma(Y_k)$-measurable, one gets
\begin{eqnarray*}
\hskip -0.25 cm \E\left(p^n_{k+1}(\chi^{n-1-k}(Q,x),Y_{k+1})\,|\F_k\right)&-&\E \left (\widehat  p^n_{k+1}(\chi^{n-1-k}(Q,x),\widehat Y_{k+1})\,|\, \widehat Y_k \right) \\
&=&\Theta_k (p^n_{k+1}(\chi^{n-1-k}(Q,x),.))(Y_k)-\E \left (\Theta_k  (p^n_{k+1}(\chi^{n-1-k}(Q,x),.)(Y_k)\,|\, \widehat Y_k\right)
\\
&&+\E\left(p^n_{k+1}(\chi^{n-1-k}(Q,x),Y_{k+1})-\widehat p^n_{k+1}(\chi^{n-1-k}(Q,x),\widehat Y_{k+1})\,|\, \widehat Y_k\right)\\
&=& \Theta_k (p^n_{k+1}(\chi^{n-1-k}(Q,x),.))(Y_k)-\Theta_k  (p^n_{k+1}(\chi^{n-1-k}(Q,x),.)(\widehat Y_k)\\
&&+ \E \left (\Theta_k  (p^n_{k+1}(\chi^{n-1-k}(Q,x),.)(\widehat Y_k) -\Theta_k  (p^n_{k+1}(\chi^{n-1-k}(Q,x),.)(Y_k) \,|\, \widehat Y_k \right)\\
&&+\E \left (p^n_{k+1}(\chi^{n-1-k}(Q,x),Y_{k+1}) - \widehat p^n_{k+1}(\chi^{n-1-k}(Q,x),\widehat Y_{k+1})\,|\, \widehat Y_k\right).
\end{eqnarray*}

Consequently, still using the elementary inequality $|\sup_{x\in X}a_x-\sup_{x\in X}b_x| \le \sup_{x\in I} |a_x-b_x|$ for any index set  $X$
and, for every $x\!\in \partial I^{n-1-k}$, that
$$
\chi^{n-k-1}(T_\N^+(n-k),x)\subset T_\N^+(n-k-1)
$$
 (see the proof of Theorem~\ref{main}$(b)$), one has
\begin{eqnarray*}
\sup_{Q\in T_\N^+(n-k)}|p^n_k(Q,Y_k)-\widehat p_k(Q,\widehat Y_k)|&\le&|v_k(Y_k)-v_k(\widehat Y_k)|\\
&&+ \sup_{Q'\in T_\N^+(n-k-1)}| \Theta_k \left(p^n_{k+1}(Q',.)\right)(Y_k)-\Theta_k  (p^n_{k+1}Q',.)(\widehat Y_k)|\\
&&+ \E \left (\sup_{Q'\in T_\N^+(n-k-1)}|\Theta_k  (p^n_{k+1}(Q',.)(\widehat Y_k) -\Theta_k  (p^n_{k+1}(Q',.)(Y_k)| \,|\, \widehat Y_k \right)\\
&&+\E \left (\sup_{Q'\in T_\N^+(n-k-1)}|p^n_{k+1}(Q',Y_{k+1}) - \widehat p^n_{k+1}(Q',\widehat Y_{k+1})|\,|\, \widehat Y_k\right).
\end{eqnarray*}

Temporarily set for convenience, $\Delta^{n,p}_k:= \|\sup_{Q\in T_\N^+(n-k)}|p^n_k(Q,Y_k)-\widehat p_k(Q,\widehat Y_k)|\|_{_p}$. One derives that for every $k=0,\ldots, n-1$,
\[
\Delta^{n,p}_k\le ([v_k]_{\rm Lip}+2[\Theta_k]_{\rm Lip}[p^n_{k+1}]_{{\rm Lip},y}) \|Y_k-\widehat Y_k\|_{_p} + \Delta^{n,p}_{k+1}.
\]
Furthermore, $\Delta^{n,p}_{n-1}\le [v_{n-1}]_{\rm Lip}\|Y_{n-1}-\widehat Y_{n-1}\|_{_p}$. The result follows by induction.$\cqfd$

\subsection{Optimal quantization}
\paragraph{Theoretical background}In this section, we provide a few basic elements about optimal quantization in order to
give some error bounds for the premium of the swing option.  We refer to~\cite{GRLU} for more details about theoretical aspects
and to~\cite{PAPR} for the algorithmic aspects numerical applications.

Let $p\!\in [1,+\infty)$. let $Y\!\in L^p(\Omega, {\cal A}, \P)$ be an  $\R^d$-valued random vector  and let    $N\ge 1$
be a given grid size. The best $L^p$-approximation of $Y$ by a random vector taking its values in a given grid $\Gamma$ of size
(at most) $N$ is given by a Voronoi quantizer $\widehat Y^\Gamma$ which induces an
$L^p(\P)$- mean quantization error
\[
e_{N,p}(Y,\Gamma) = \|Y-\widehat Y^\Gamma\|_{_p}=
\left(\E\min_{y\in \Gamma}|Y-y^i|^p \right)^{\frac 1p}.
\]
It has been shown independently by several authors (in various finite and infinite dimensional frameworks) that when the
grids $\Gamma$ runs over all the subsets of $\R^d$ of size at most $N$, that $e_{N,p}(Y,\Gamma)$ reaches a minimum denoted
$e_{N,p}(Y)$ (see~$e.g.$~\cite{GRLU} or~\cite{PAG}) $i.e.$ the minimization problem
\[
e_{N,p}(Y)= \min\left\{e_{N,p}(Y,\Gamma), \; \Gamma\subset \R^d,\; \mbox{\rm card}(\Gamma)\le N  \right\}
\]
has at least a solution temporarily denoted $\Gamma^{(N,*)}$. Several algorithms have been designed to
compute some optimal or close to optimality quantizers, especially in the quadratic case $p=2$. They all rely on the stationarity
property
satisfied by optimal quantizers. In the quadratic case, a grid $\Gamma$ is stationary
\[
\widehat Y^{\Gamma}= \E\left(Y\,|\, \widehat Y^{\Gamma}\right)
\]
This follows from some differentiability property of the $L^p$-distortion. For a formula in the general case we
refer to~\cite{GRLUPA3}. In $1$-dimension, a regular Newton-Raphson zero search procedure turns out to be quite efficient. In higher dimension
(at least when $d\ge 3$ or $4$) only stochastic procedures can be implemented like the $CLVQ$ (a stochastic gradient descent, see~\cite{PAG} or
the Lloyd~I procedure (a randomized fixed point procedure, see~\cite{IEEE}). For more details and result we refer to~\cite{PAPR}.

As a result of these methods, some optimized grids of the (centered) normal distribution ${\cal N}(0;I_d)$ are available
on line at the URL
\[
\mbox{\tt www.quantize.maths-fi.com}
\]
for dimensions $d=1,\ldots,10$ and sizes from $N=2$ up to $5\, 000$.

\ms It is clear by considering a sequence of grids $\Gamma^{(N)}:=\{r^1,\ldots,r^{N}\}$  where $(r^n)_{n\ge1}$ is an everywhere
dense sequence in $\R^d$ that  $e_{N,p}(Y)$ decreases to $0$ as $N\to \infty$.

The rate of convergence of this sequence is ruled by the so-called Zador Theorem (see~\cite{ZAD1} for a first statement of the result, until
the first rigorous proof in~\cite{GRLU}).

\begin{Thm} \label{Zad}(Zador, see~\cite{GRLU}) $(a)$ Let $Y\!\in L^{p+\eta}(\P)$, $p\ge 1$, $\eta>0$, such that $\P_{_Y}(du)= \varphi(du)du\stackrel{\perp}{+}\nu(du)$. Then
\[
\lim_N N^{\frac 1d} e_{N,p}(Y)= \widetilde J_{p,d} \left(\int_{\R^d} \varphi^{\frac{d}{p+d}}(u)\,du\right)^{\frac{1}{p}+\frac{1}{d}}.
\]
$(b)$ Non asymptotic estimate (see $e.g.$~\cite{LUPA4}):  Let $p\ge 1$, $\eta >0$. There exists a real constant
$C_{d,p, \eta}>0$ and an integer $N_{d,p,\eta}\ge 1$ such that or any $\R^d$-valued random vector $Y$, for $N\ge N_{d,p,\eta}$,
\[
e_{N,p}(Y)\le C_{d,p, \eta}\|Y\|_{p+\eta}N^{-\frac 1d}.
\]
\end{Thm}


\paragraph{Rate of convergence of the quantization pricing method} Now we are in position to apply the above results to provide an error bound for the pricing of swing options by optimal quantization: assume there is a real exponent $p\!\in [1,+\infty)$ such that the ($d$-dimensional) Markov structure process $(Y_k)_{0\le k\le n-1}$ satisfies
\[
\max_{0\le k\le n-1}|Y_k|\!\in L^{p+\eta}(\P), \;\eta>0
\]
At each time $k\!\in\{0,\ldots,n-1\}$, we implement a (quadratic) optimal quantization grid $\Gamma^{\bar N}$ of  $Y_k$ with constant  size
$\bar N$ . Then the  general error bound result~(\ref{ratequantswing}) combined with  Theorem~\ref{Zad}$(b)$ says that, if   $\bar N \ge
N_{d,p,\eta}$,
\begin{eqnarray*}
\left\|\sup_{Q\in T^+(n)}|P^n_0(Q) -  \widehat  P^n_0(Q)|\right\|_{_p}&\le& C \frac{n}{\bar N^{\frac 1d}}
\end{eqnarray*}
where as the complexity of the procedure is bounded by $n(n+1)\bar N^2$ (up to a constant).

In fact this error bound turns out to be conservative and several numerical experiments, as those presented below,   suggest that in fact the true rate
(for a fixed number $n$ of purchase instants)  behaves like $O(\bar N^{-\frac 2d})$.

Another approach could be to minimize the complexity of the procedure by considering (optimal) grids with variable sizes $N_k$ satisfying
$\sum_{k=0}^{n-1} N_k= n\,\bar N$. We refer to~\cite{BABOPA1} for further results in that direction. However, numerical experiments were
  carried out with constant size grids for both programming convenience and memory saving.

\subsection{A numerical illustration}

We considered a two factor continuous model for the price of future contracts which
leads to the following dynamics for the spot price
\[
S_t = F_{0,t} \exp{\left(\sigma_1\int_0^te^{-\alpha_1(t-s)}dW^1_s +
\sigma_2\int_0^t e^{-\alpha_2(t-s)}dW^2_s -
\frac{1}{2}\Lambda_t\right)},\qquad t\!\in[0,T],
\]
where $W^1$ and $W^2$ are two standard Brownian motions with
correlation coefficient $\rho$ and
\[\Lambda_t = \frac{\sigma_1^2}{2\alpha_1}\left(1-e^{-2\alpha_1
t}\right)+\frac{\sigma_2^2}{2\alpha_2}\left(1-e^{-2\alpha_2
t}\right)+\frac{2\rho\sigma_1\sigma_2}{\alpha_1+\alpha_2}\left(1-e^{-(\alpha_1+\alpha_2)
t}\right)\]

Then, we consider a (daily) discretization of  the Gaussian process
$\log(S_t/F_{0,t})$ at times $\frac{kT}{n}$ $(T=1$, $n=365$). The
sequence $(\log(S_{t_k}/F_{0,t_{k}}))_{0\le k\le n-1}$ is clearly
not Markov. However, adding an  appropriate auxiliary processes, one
can build a higher dimensional (homogenous) Markov process
$(Y_k)_{0\le k \le n-1}$ whom $\log(S_{t_k}/F_{0,t_{k}})$ is a
linear combination. This calls upon classical methods coming from
time series analysis. Then a fast quantization method has been
developed to make makes possible a parallel implementation of the
quantized probability transitions of $Y=(Y_k)_{0\le k\le n-1}$. For
further details about this model and the way it can be quantized, we
refer to~\cite{BABOPA1}. In~\cite{BABOPA1},  the optimal
quantization method described above is extensively tested from a
numerical viewpoint (rates of convergence, needed memory, swapping
effect, etc). Its performances are compared those of the
Longstaff-Schwartz approach introduced in~\cite{GobetDF}. This
comparison emphasizes the accuracy and the velocity of our approach,
even if only  one contract is to be computed and  the computation of
the probability transitions is included in the computation time of
the quantization method. Furthermore, it seems that it needs
significantly  less memory capacity when implemented on our tested
model.

We simply reproduced below a complete graph of the function
$Q:=(Q_{\min}, Q_{\max})\mapsto P^n_0(Q)$  when $Q$ runs over the whole set of admissible global constraints $T^+(n)$. The parameters were settled at the following values
\[
n=30,\; \alpha_1= 0.21,\; \alpha_2=5.4,\; \sigma_1=36\%,\; \;
\sigma_2=111\%,\;\rho=-0.11\]

The graph of the premium function $Q\mapsto P^n_0(Q)$ defined on $T^+(n)$  is depicted in Figure~1.


\begin{figure}
\centerline{\includegraphics[angle=0,width=12.5cm]{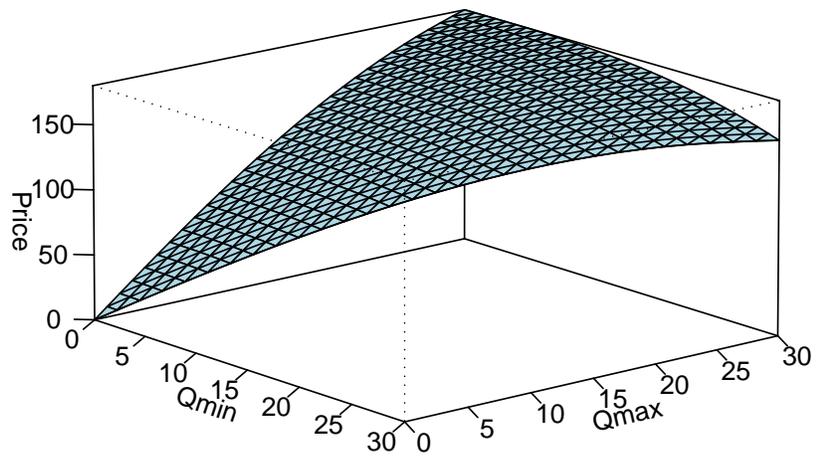}}
\caption{{\em The mapping $Q\mapsto \widehat P^n_0(Q)$ affinely
interpolated from integral-valued global constraints.}}
\label{fig:2}
\end{figure}

\bigskip
\small
\noindent {\sc Acknowledgement:} We thank Anne-Laure Bronstein for helpful comments.

\end{document}